\documentclass[11pt,reqno]{amsart}

\usepackage{notation2}

\usepackage{adjustbox}

\title{On Nash-Williams' Theorem regarding sequences with finite range}
\author{Fedor Pakhomov and Giovanni Soldà}

\address{Fedor Pakhomov, Department of Mathematics: Analysis, Logic and Discrete Mathematics, Ghent University, Krijgslaan 281 S8, 9000 Ghent, Belgium, and\newline
Steklov Mathematical Institute of the Russian Academy of Sciences. Ulitsa
Gubkina 8, Moscow 117966, Russia}
\email{fedor.pakhomov@ugent.be}

\address{Giovanni Sold\`a, Department of Mathematics: Analysis, Logic and Discrete Mathematics, Ghent University, Krijgslaan 281 S8, 9000 Ghent, Belgium}
\email{giovanni.a.solda@gmail.com}

\thanks{The work of the authors has been funded by the FWO grant G0F8421N}

\begin{document}

\begin{abstract}
The famous theorem of Higman states that for any well-quasi-order (wqo) $Q$ the embeddability order on finite sequences over $Q$ is also wqo. In his celebrated 1965 paper, Nash-Williams established that the same conclusion holds even for all the transfinite sequences with finite range, thus proving a far reaching generalization of Higman's theorem. 

In the present paper we show that Nash-Williams' Theorem is provable in the system $\mathsf{ATR}_0$ of second-order arithmetic, thus solving an open problem by Antonio Montalbán and proving the reverse-mathematical equivalence of Nash-Williams' Theorem and $\mathsf{ATR}_0$. In order to accomplish this, we establish equivalent characterization of transfinite Higman's order and an order on the cumulative hierarchy with urelements from the starting wqo $Q$, and find some new connection that can be of purely order-theoretic interest. Moreover, in this paper we present a new setup that allows us to develop the theory of $\alpha$-wqo's in a way that is formalizable within primitive-recursive set theory with urelements, in a smooth and code-free fashion.
\end{abstract}

\maketitle

\section{Introduction}\label{sec:introduction}

A \emph{well quasi-order} is a quasi-order without infinite descending sequences and infinite antichains. Equivalently, an order $Q$ is well quasi-order if it does not admit any \emph{bad sequence}, i.e.\ a sequence $f:\omega\to Q$ such that $i<j$ entails $f(i)\not\leq_Q f(j)$.
This notion is arguably one of the most natural and fundamental in mathematics: the concept has been frequently reintroduced and rediscovered in different contexts, to the point of prompting Kruskal to write \cite{wqo-freq-disc-kruskal} ``so that future work will not needlessly repeat what is known''. 

At present, well quasi-orders are at the core of a rich theory with many applications. An important part of this theory is constituted by results that give closure properties for well quasi-orders: among them, of capital importance are Higman's Lemma, first proved in \cite{hig-theo-higman} (and later famously generalized by Kruskal in \cite{kruskal-theo-kruskal}), which states that if $Q$ is well quasi-order, then the set of finite sequences over $Q$ is also a well quasi-order, if we order it with the embeddability relation (we refer to \Cref{def:all-seq} for a precise formulation of this relation). It was soon noticed, anyway, that this result is in a way optimal: it is not the case in general that sequences of length $\omega$ over a well quasi-order $Q$ form a well quasi-order, as witnessed for instance by Rado's counterexample (\cite{counterex-wqo-rado}). It was the pioneering work of Crispin Nash-Williams that gave a more comprehensive view of the relationship between the properties of a quasi-order $Q$ and those of the transfinite sequences on $Q$. On the one hand, in \cite{wqo-transf-seq-nash-williams} (building on work of Kruskal, Erdős and Rado) he showed that if $Q$ is well quasi-ordered, then the class of transfinite sequences \emph{with finite range} is well quasi-ordered by embeddability.  And on the other, in \cite{nwt-nash-williams} he gave a classification of the quasi-orders $Q$ such that the whole class of transfinite sequences over $Q$ forms a well quasi-order: 
such $Q$'s are exactly the \emph{better quasi-orders}, a class of orders first introduced in \cite{wqo-inf-trees-nash-williams}. 

Since their inception, better quasi-orders have been a very important tool of infinitary combinatorics: they enjoy much stronger closure properties than well quasi-orders, which makes it somewhat easier to work with them. Famously, this was exploited in the proof given by Laver in \cite{fraisee-conj-laver} of the fact that countable scattered linear orders form a better quasi-order (a strengthening of Fraïssé's conjecture). This comes at a price: the original definition of a better quasi-order is rather complicated (and it is then usually difficult to verify that a certain order is a better quasi-order). We postpone the official definition to \Cref{def:wqo-bqo}, but we give an intuitive approximation of it here: just as well quasi-orders were orders not admitting any bad sequence, better quasi-orders do not admit any \emph{bad array}, which are functions having as domains certain subsets of $[\omega]^{<\omega}$ called \emph{fronts}. Fronts intuitively correspond to minimal clopen covers of Baire space (this was exploited by Simpson to give a ``simpler" description of better quasi-orders in \cite{bqo-fraisse-simpson}), and can thus be seen as well-founded trees: so to every front $F$ we can assign a rank $\alpha$, namely the rank of the associated well-founded tree. This perspective makes it possible to attach a ``measure of better quasi-orderedness" to an order $Q$: we can ask what is the maximal $\alpha$ such that no front of rank at most $\alpha$ is the domain of a bad array with range $Q$, and we will call \emph{$\beta$-well quasi-orders} the quasi-orders for which this $\alpha$ is at least $\beta$. This measure has important order-theoretic properties, as we shall see later on.

From the point of view of mathematical logic, the theory of well and better quasi-orders is quite interesting: the survey \cite{survey-wqo-bqo-2020-marcone} gives an account of this. Notably, several questions about the logical strength of theorems concerning well and better quasi-orders remain wide open. Our paper is a contribution to this line of research: in particular, we show that Nash-Williams' result that if $Q$ is well quasi-ordered, then so is the class of transfinite sequences over $Q$ with finite range, is provable in $\mathsf{ATR}_0$. Since it follows from Shore's work in \cite{fraisse-conj-shore} (together with the results in \cite{weak-well-fraisse-freund-manca}, where a small gap in Shore's argument is identified) that this theorem implies $\atr$ over $\rca$, this determines its reverse mathematical strength, thereby answering Question $23$ of \cite{open-questions-montalban}.

To achieve this we develop a framework for the formalization of the theory of $\alpha$-well quasi-orders. Namely, in \Cref{sec:tr-seq}, we consider a number of different characterizations of the property of a quasi-order $Q$ to be an $\alpha$-well quasi-order (the following is a simplification of the content of \Cref{theo:big-theo-1} and \Cref{theo:big-theo-2}):
\begin{itemize}
    \item $Q$ has no bad arrays of the ranks $<\alpha$;
    \item the level  $\dot V_\alpha(Q)$ of cumulative hierarchy with urelements from the order $Q$ is well-founded according to certain natural extension of $Q$ (where the dot indicates that we restrict the cumulative hierarchy to urelements and sets that do not contain $\emptyset$ in their transitive closure);
    \item the embeddability order on transfinite sequences of the length $<\omega^{1+\alpha}$ is well-founded.
\end{itemize}
Furthermore we establish that under any of the conditions above, the quasi-order on the indecomposable transfinite sequences of the length $<\omega^{1+\alpha}$ and the quasi-order $\dot V_\alpha(Q)$ are isomorphic after the factorization by their respective internal equivalence relations. 

We point out that the observation that there exists a tight relationship between the degree of better quasi-orderedness of an order $Q$ and the order-theoretic properties of the iterated powerset of $Q$was given already in \cite{wq-trees-nash-williams} by Nash-Williams, and expanded in \cite{fraisee-conj-laver} (in addition to the previous references, we also mention \cite{phd-pequignot} for a very clear explanation of this relationship). The connection between degree of better quasi-orderedness and maximal $\gamma$ such that $Q^{<\gamma}$ is well-founded was also noticed before, and we refer to \cite{theory-relations-fraisse} and \cite{foundations-bqo-marcone} for a treatment of this matter.

In \Cref{sec:tr-seq-fin-ran}, we find that the relationship between iterated powerset and transfinite sequences can be fruitfully applied to the case of sequences with finite range (since they could be considered to be sequences over a finite order), and allows us to reduce Nash-Williams' result on sequences of finite range to the statement that for a well-quasi orders $Q$ there are no bad sequences in $\dot V_\alpha(Q)$ restricted to elements with finite support. Using the general connection between bad sequences in $\dot V_\alpha(Q)$ and bad arrays, any bad sequence as in the previous sentence could be transformed into a bad array over $Q$. At this point we are able to finish the proof of Nash-Williams' Theorem by showing that the specific arrays as above could be transformed into bad sequences over $Q$.

We develop the aforementioned theory in the way that is directly formalizable in a variant $\prsou(\mathsf{enu})$ of primitive recursive set theory. Namely the theory $\prsou(\mathsf{enu})$ is a set theory with urelements and a global enumerating function, i.e. the function picking a countable enumeration for every non-empty set (in particular this implies that the theory proves all sets to be at most countable). We transfer this formalization to $\mathsf{ATR}_0$ by a two-fold construction. First we directly and naturally interpret in $\atr$ the variant of the theory above, where the global enumerating function is replaced with the axiom of countability\footnote{The interpretation is basically the urelement version of the known interpretation of $\mathsf{ATR}_0^{\text{set}}$ in $\atr$ \cite{simpson1982set,book-simpson}. We note that the main difference between $\mathsf{ATR}_0^{\text{set}}$ and primitive recursive set theory is that the former has axiom $\beta$ (Mostowski collapsing theorem considered as an axiom) and the latter does not.}. And secondly, via a forcing construction we show that any argument involving the global enumerating function but not using it in its conclusion could be transformed into an argument that uses just the axiom of countability (this is done in \Cref{app:forcing-enu}).

\section{Preliminaries}

\subsection{The theory $\prsou$}\label{sec:prsou}

Within the present paper we will be mostly working in primitive recursive set theory with urelements $\prsou$. For a more gentle presentation and a proof of several results that we simply state, we refer to the first chapter of \cite{phd-freund}\footnote{Technically, what Anton Freund was dealing with was the theory $\prso$, but it is rather clear that the presence of urelements does not change much}. The language of this theory extends the language of set theory by additional constants $\emptyset,\omega$, and function symbols for primitive-recursive set functions:
\begin{enumerate}
\item $C(x,y,z,w)=x$ if $z\in w$ and $C(x,y,z,w)=y$ otherwise;
\item $A(x,y)=x\cup\{y\}$;
\item $R[t(x,y,\vec{z})](y,\vec{z})=t(\bigcup\{R[t(x,y,\vec{z})](y',\vec{z})\mid y'\in y\},y,\vec{z})$, where $t(x,y,\vec{z})$ is a primitive-recursive set term.
\item $U(x)$ is the set of all urelements that are elements of $x$.
\end{enumerate}
The bounded quantifiers are interpreted in the standard way, expression $\mathsf{U}(x)$ reads as $x$ is an urelement and is an abbreviation for $x\in U(A(\emptyset,x))$.
The axioms are
\begin{enumerate}
\item $x\not \in \emptyset$
\item $y\in U(x)\to y\in x$\
\item $z\in A(x,y)\mathrel{\leftrightarrow} z\in x\lor z=y$
\item $\lnot \mathsf{U}(x)\land \lnot\mathsf{U}(y)\land \forall z\in x(z\in y) \land \forall z\in y (z\in x)$ (Extensionality)
\item $\mathsf{U}(x) \to y\not\in x$
\item $\lnot\mathsf{U}(\emptyset)$
\item $\forall y (\forall z\in y\;\varphi(z)\to\varphi(y))\to\varphi(x)$, for $\Delta_0$ formulas $\varphi$
\item $x\in\omega\mathrel{\leftrightarrow} x=\emptyset \lor \exists y\in x\;x=A(y,y)$
\item $z\in w\to C(x,y,z,w)=x$
\item $z\not \in w \to C(x,y,z,w)=y$
\item $\emptyset\in \omega\land (x\in \omega \to A(x,x)\in \omega)$
\item $\forall x\in \omega (x=\emptyset \lor \exists y\in \omega(x=A(y,y)))$
\item $\exists w\Big (\big ( \forall v\in w\exists u\in y\; v=R[t(x,y,\vec{z})](u,\vec{z})\big)\\\land\;\big( \forall u\in y\exists v\in w\;v=R[t(x,y,\vec{z})](u,\vec{z})\big )\\\land R[t(x,y,\vec{z})](y,\vec{z})=t(w,y,\vec{z})\Big )$
\end{enumerate}

The constructions carried out in $\prsou$ will in fact also work in $\mathsf{ATR}_0$. This is because we could naturally interpret $\prsou$ in $\mathsf{ATR}_0$ with  $U=\mathbb{N}$.

We will be more precise about this interpretation in \Cref{app:prosu-in-atr}, but we nevertheless give an intuitive description of how this interpretation works. First of all, we give a slight modification of the standard interpretation of $\mathsf{ATR}_0^{\mathsf{set}}$ in $\mathsf{ATR}_0$: sets are interpreted as pairs $((A,\in^A),a)$, where $(A,\in^A)$ is a well-founded extensional binary relations and $a\in A$. $((A,\in^A),a)$ is equal as a set to $((B,\in^B),b)$ if there exists a well-founded extensional amalgamation $(C,\in^C)$ of $(A,\in^C)$ and $(B,\in^B)$ by initial embeddings $f_A\colon (A,\in^A)\to (C,\in^C)$ and $f_B\colon (B,\in^B)\to (C,\in^C)$ such that $f_A(a)=f_B(b)$. Membership of $((A,\in^A),a)$ in $((B,\in^B),b)$ is interpreted analogously, but now we demand that $f_A(a)\in f_B(b)$. Although this might not look familiar, we observe that we essentially just recast the argument of \cite[Section VII.3]{book-simpson} in a slightly different fashion: if $((A,\in^A),a)$ is meant to represent the set $u$, then one can see $A$ as the suitable tree associated to $u$, $\in^A$ as the opposite of the successor function on $T$, and $a$ as the root of $T$.

Now we modify this interpretation for the case when we have $\mathbb{N}$ as the set of urelements. I.e.\ now, instead of well-founded extensional binary relations $(A,\in^A)$, we consider binary relations with pairwise distinct constants $\mathsf{nmb}_n$, for $n\in\mathbb{N}$ that are extensional for elements that are not these constants. 

We were, however, not able to find a satisfactory reference for the interpretation of the primitive recursive set functions: this is not particularly difficult, and seems to be dealt with in different places as a matter of course. It nevertheless seems to require some thought: we give such an interpretation in Appendix \ref{app:prosu-in-atr}.

We notice that the interpretation above satisfies the axiom ``every set is countable'', i.e. that for every set $x$ there is a function $f$ with $\mathsf{dom}(f)=\omega$ and $\mathsf{ran}(f)\supseteq x$.

In the following, it will sometimes be practical to work in the theory with the global enumerating function.

\begin{Definition}$\prsou(\mathsf{enu})$ is the extension of $\prsou$ with one additional binary function $\mathsf{enu}(x,y)$ that satisfies the axiom
  $$x\ne \emptyset\land \neg \mathsf{U}(x) \to \forall x'\in x\exists y\in\omega(\mathsf{enu}(x,y)=x')\land \forall y\in \omega(\mathsf{enu}(x,y)\in x).$$
\end{Definition}

Under the axiom of countability, there is very little difference between $\prsou$ and $\prsou(\mathsf{enu})$: the following result is proved in \Cref{app:forcing-enu}.

\begin{Lemma}$\prsou(\mathsf{enu})$ is a conservative extension of $\prsou+\text{``every set is countable''}$.\end{Lemma}


\subsection{Basic definitions of wqo- and bqo-theory}

We now give, over the theory $\prsou$,
some important definition that will be fundamental for the
following. Several of them are standard, but we include
them nevertheless: we do this to fix notation and to spell
out some peculiarity that our choice of base theory entails.

\begin{Definition}
($\prsou$)
  A pair $(Q,\leq_Q)$ is a \emph{quasi-order} (henceforth \emph{qo}) if $Q$ is
  a primitive recursive class and $\leq_Q$ is a primitive recursive binary transitive and reflexive relation on it.
\end{Definition}

In most of the following, $Q$ will be  the class of urelements. We will often abuse
notation and refer to $Q$ as a qo without specifying the relation $\leq_Q$,
as customary in the literature.

\begin{Definition}
    ($\prsou$)
  For a set of natural numbers $A$, the set $[A]^{<\omega}$ is the set of all strictly ascending finite sequences built from the elements of $A$. The order $\sqsubseteq$ on $[\omega]^{<\omega}$ is the extension order.
\end{Definition}

There is a natural bijection between $[A]^{<\omega}$ and $\mathsf{P}_{\mathsf{fin}}(A)=\{x\subseteq A\mid x \text{ is finite}\}$ that maps each ascending sequence $(a_0,\ldots,a_{n-1})$ to the finite set $\{a_0,\ldots,a_{n-1}\}$. Thus we will often abuse notations and identify ascending sequences with the corresponding finite sets.

We now reach the first point where our choice of base theory matters:
the definition of fronts. The problem for us is the existence
of \emph{ranks} on them: unlike $\atr$, $\prso$ is not strong enough
to prove the existence of ranks for arbitrary fronts. We take care of this
by requiring that every front already comes equipped with a
suitable rank.

\begin{Definition}\label{def:front}
($\prsou$)
  A \emph{front} is a pair $(F,\mathsf{rk}_F)$ 
  where $F\subseteq [\omega]^{<\omega}$ is an infinite $\sqsubseteq$-antichain
  such that, for every $X\in [\bigcup F]^{\omega}$, there is $\sigma\in F$
  such that $\sigma\sqsubset X$, 
  and $\mathsf{rk}_F$ is a function from
  $$T_F=\{\sigma \in \bigl[ \bigcup F \bigr]^{<\omega}\mid
  \exists \tau \in F( \sigma\sqsubset \tau)\}$$
  to ordinals such that $\mathsf{rk}_F(\sigma)>\mathsf{rk}_F(\tau)$  for
  $\sigma\sqsubset\tau$ from $T_F$. 

  We say that the \emph{rank of $(F,\mathsf{rk}_F)$ is $\le \alpha$} if
  $\mathsf{rk}_T(\emptyset)\le\alpha$, and that it is $<\alpha$ if
  $\mathsf{rk}_F(\emptyset)<\alpha$.
\end{Definition}

We remark that in the definition above by ``$\sigma\sqsubset \tau$" we mean ``$\sigma\subseteq\tau$ and $\sigma\neq\tau$".

We will often abuse the notations and identify a front $(F,\mathsf{rk}_F)$  with $F$ or in other words assume that when we considering $F$, we are also implicitly fixing a rank function $\mathsf{rk}_F$.

In the following, we will need to ``restrict'' fronts on sets smaller than $\bigcup F$.
By this we mean the following: given a front $F$ and and an infinite
$H\subseteq \bigcup F$ we consider the front
 $$F\rst_H=\{\sigma\in F: \sigma\subseteq H\},$$ where $\mathsf{rk}_{F\rst_H}=\mathsf{rk}_F\rst_{F\rst_H}$. Note that $\prsou$ proves that under the considered conditions $F\rst_H$ will be a front.

\begin{Definition}
($\prsou$)
  For a non-empty $\sigma\in[\omega]^{<\omega}$ we put
  $$\sigma^-=\sigma\setminus \min \sigma.$$
  For non-empty $\sigma,\tau\in[\omega]^{<\omega}$ we put
  $$\sigma\triangleleft \tau  \stackrel{\mbox{\scriptsize\textrm{def}}}{\iff}\min(\sigma)<\min(\tau) \text{ and }\sigma^- \text{ is }\sqsubseteq\text{ comparable with }\tau$$
\end{Definition}

\begin{Definition}
($\prsou$)
  An \emph{array} is a function $f\colon F\to P$ from a front $F$ to a quasi-order $P$. We say that $f$ has \emph{rank} $\le\alpha$ if $F$ has rank $\le \alpha$. We say that $f$ is \emph{bad} if $f(\sigma)\not\le_P f(\tau)$ holds for every $\sigma,\tau\in F$, $\sigma\triangleleft \tau$. We say that $f$ is \emph{good} if it is not bad.
\end{Definition}

\begin{Definition}\label{def:wqo-bqo}
($\prsou$)
  We say that a quasi-order $Q$ is an \emph{$\alpha$-well quasi-order} (henceforth \emph{\emph{$\alpha$}-wqo}) if it is well-founded and there are no bad arrays $f\colon F\to Q$ of the rank $<\alpha$.

  We say that that $Q$ is a \emph{better quasi-order} (henceforth \emph{bqo}) if it is $\alpha$-wqo for any $\alpha$.
\end{Definition}

Note that fronts $F$ of the rank $\le 0$ are precisely infinite sets of singletones.
This in particular allows us to recover the classical concept of wqo as a particular
case of the above: namely, $Q$ is wqo if and only if it is $1$-wqo. As a further
particular case, notice  that $0$-wqo's are precisely the well-founded qo's.

\begin{Definition}\label{def:all-seq}
($\prsou$)
  A \emph{transfinite sequence} $u$ over a qo $Q$ of \emph{length} $|u|\in\mathsf{On}$, $|u|>0$ is a function with domain $|u|$ and range contained in $Q$. We identify $q\in Q$ with the sequence $u$ of length $|u|=1$ and $u(0)=q$. Naturally we define sums $u+v$ of sequences $u,v$ as follows: $|u+v|=|u|+|v|$ and $(u+v)(\alpha)=u(\alpha)$, for $\alpha<|u|$ and $(u+v)(|u|+\alpha)=v(\alpha)$, for $\alpha<|v|$. For an indexed family of transfinite sequences $(u_\beta)_{\beta<\alpha}$ the sequence $\sum_{\beta<\alpha} u_\beta$ is defined in the obvious way.

  For two transfinite sequences $u,v$ a \emph{(weak) embedding} $f\colon u\hookrightarrow v$ ($f\colon u\hookrightarrow_w v$) is a strictly (non-strictly) increasing function $f\colon |u|\to |v|$ such that $u(\beta)\le_Q v(f(\beta))$.  We denote as $\seq \alpha Q$ ($\seqalt \alpha Q$) the (weak) embeddability order on sequences of the length $<\alpha$, and we write $u\preceq v$ ($u\preceq^* v$) if there exist a
  (weak) embedding $f\colon u\hookrightarrow v$ ($f\colon u\hookrightarrow_w v$).

  Given ordinals $\alpha<\beta<|u|$, we denote by $u\rst_{[\alpha,\beta)}$ the
  restriction of $u$ on the semiclosed interval $[\alpha,\beta)$, namely
  the sequence $v$ whose domain is $\{\gamma\in\mathsf{On}:\alpha+\gamma<\beta\}$
  and such that $v(\gamma)=u(\alpha+\gamma)$. Any $u\rst_{[\alpha,\beta)}$ is called a \emph{subsegment} of $u$. Subsegments of the form $u\rst_{[\alpha,|u|)}$ are called \emph{tails}.

  A sequence $u$ is called \emph{(weakly) indecomposable} if for any $\alpha<|u|$, $u$ (weakly) embedds into $u\rst_{[\alpha,|u|)}$, i.e.\ if $u$ (weakly) embeds into all of its tails. The suborder of $\seq \alpha Q$ ($\seqalt \alpha Q$) on (weak) indecomposable sequences is denoted as $\indseq \alpha Q$ ($\indseqalt \alpha Q$).
\end{Definition}

Note that the (weak) embeddability relation between seqences in fact is primitive set recursive. This is because if there is an embedding $f\colon v\hookrightarrow u$, then the following recursively defined function $f'\colon |v|\to |u|$ is an embedding:
$$f'(\alpha)=\mathsf{min} \{\beta \mid \forall \alpha'<\alpha (\beta>f'(\alpha'))\text{ and } v(\beta)\le_Q u(\beta)\},$$
where, in the case of weak embeddings, the condition $\beta>f'(\alpha')$ should be replaced with $\beta\ge f'(\alpha')$. This is essentially \cite[Lemma 3.4]{nash-williams-marcone}. 

\subsection{Two important tools: Higman's Lemma and Clopen Ramsey's Theorem}
Here we recall two important instruments of wqo- and bqo-theory, namely Higman's
Lemma and Clopen Ramsey's Theorem, and some of their consequences. We focus on
the fact that they are indeed available in the base theory $\prsou$: whereas
this is very much not surprising in the case of Higman's Lemma, the proof of
Clopen Ramsey's Theorem relies on our careful definition of fronts.
We will not give full proofs, but sketches highlighting some
of the differences between the arguments in second-order arithmetic
and $\prsou$.

\begin{Theorem}[Higman's Lemma]
($\prsou$)
  If $Q$ is a wqo, then so is $\seq \omega Q$.
\end{Theorem}
\begin{proof}
We first notice that we can restrict to the case of countable $Q$: if $Q^{<\omega}$ admits a bad sequence $f$, then we can consider $\{ q\in Q: \exists i,j (q=f(i)(j))\}$, which is countable, in place of $Q$.
  Then, the proof proceeds in the standard way (see, for instance,
  \cite[Theorem 3]{metamat-fraisse-clote} and
  \cite[Sublemma 4.8]{hilbert-basis-simpson}): given a wqo $Q$, we can create
  the tree $\text{Bad}(Q)$ of its bad sequences: since we can prove that
  it is well-founded, we can conclude that if we order it by the
  Kleene-Brower order we obtain a well-order, say
  $(L+1,\leq_L)$
  (we are compelled to carry the ``$+1$'' in order to adhere
  to the bibliographical sources we have given). We can then
  proceed to see the identity map $f:\text{Bad}(Q)\to (L+1,\leq_L)$
  as a reification of $Q$, and conclude that
  $\seq \omega Q$ is reified by $\omega^{\omega^{L+1}}$. Since $\prso$
  proves that $L$ is well-founded if and only if $\omega^L$ is,
  and that if a qo $Q$ has a reification by a well-order
  then $Q$ is wqo,
  we are done.
\end{proof}
In the proof above, we remark that, even if our
theory is sufficient to define
the well-order $L$, it does not prove that it
can be collapsed to an ordinal (which requires
axiom $\beta$). Nevertheless, this was of no consequence
in this setting, since we only cared about well-foundedness
of $L$ and the easily defined terms $\omega^{\omega^{L+1}}$, and we
did not need to perform recursion along these orders.

\begin{Theorem}[Clopen Ramsey's Theorem]
($\prsou$)
  Let $(F,\mathsf{rk}_F)$ be a front, and let $c:F\to 2$ be a coloring
  of it in two colors. Then, there exists an infinite $H\subseteq
  \bigcup F$ such that $F\rst_H$ is $c$-homogeneous.
\end{Theorem}
\begin{proof}
  Here, we follow the proof given in \cite{rec-clopen-ramsey-clote},
  by recursion along the lexicographic ordering of $T_F$. The point is
  that now, thanks to the existence of the rank function, we can prove
  that this order can be collapsed to an ordinal: indeed, if $\alpha$
  was the range of $\mathsf{rk}_F$, we can prove that there is an
  embedding of $(T_F,\leq_{lex})$ into $\omega^\alpha$. 
\end{proof}

The result above allows us to conclude that $\prsou$ can prove, in
essence, the same facts about bqo's as $\atr$.
Due to its importance, we explicitly spell out the following fact.

\begin{Corollary}
  ($\prsou$) If $P$ and $Q$ are $\alpha$-wqo's, then so
  is their product $P\times Q$.
\end{Corollary}
\begin{proof}
  Suppose not: then there is a bad $f:F\to P\times Q$, where
  the rank of $F$ is $<\alpha$. By the Clopen Ramsey's Theorem,
  we can find $H\subseteq \bigcup F$ such that
  either $\pi_P\circ f\rst_{F\rst_H}$ or $\pi_Q\circ f\rst_{F\rst_H}$
  is bad (where $\pi_P$ and $\pi_Q$ are the projections on $P$ and $Q$,
  respectively). Since $\mathsf{rk}_F$ is a rank on $F\rst_{H}$,
  this contradicts the assumptions on $P$ or $Q$.
\end{proof}

\section{Orders on the cumulative hierarchies over
  quasi-orders}

Let us fix a qo of urelements $(Q,\leq_Q)$.
In this Section, we introduce several interesting subclasses of $V(Q)$,
and corresponding orderings on them,
intuitively corresponding to ``iterated non-empty downward-closed sets''. As we will
see, there are several sensible options to formalize this idea.

We start with $\vdc Q$, the one that should be closer to 
the immediate intuition of what these
subclasses should be. We remark one important fact, on which we will
comment more extensively later: although we will give the definition
in $\prsou$, the formula corresponding to it will be 
$\Pi_1$, and so too complex to work with in weak theories.

In essence, we would like to proceed by recursion
along the ordinals as follows: we set $\idcs 0 Q =Q$, ordered by $\leq_Q$, 
and then we let 
$$\idcs{\alpha}{Q}=Q\cup \bigcup\limits_{\beta<\alpha}\mathcal{D}_Q(\idcs{\beta}{Q}),$$
where $\mathcal{D}_Q$ is the operation mapping a quasi-ordered set to the set of its non-empty, downward-closed subsets. Here for the intended behaviour we furthermore need to assume that $Q$ is set-sized so that every downwards closed collection in the hierarchy would be guaranteed to be set-sized. We order $\idcs{\alpha}{Q}$ by keeping the order $\le_Q$ on $Q$, comparing urelements with sets by the membership relation, and comparing sets between each other by the subset relation. 

Before we do that, we make an important remark: for reason that will become apparent later, \emph{we work with a notion of rank that differs from what seems to be the standard one,} i.e.\ the one defined in, say, \cite{book-barwise}. In our case, the rank only has \emph{the class of sets} as its domain, and for every set $x$ and ordinal $\alpha$ we have
\[
    \mathsf{rk}(x)\leq\alpha \Leftrightarrow \forall y\in x(y\not\in Q \rightarrow \mathsf{rk}(y)<\alpha).
\]
So, for instance, $\mathsf{rk}(\emptyset)=\mathsf{rk}(\{q\})=0$ for every $q\in Q$. We use this definition of rank to build the classes $V_\alpha(Q)$ in the usual way, namely
\[
    V_\alpha(Q) = Q \cup \bigcup_{\beta<\alpha} \mathcal{P}(V_\beta(Q)),
\]
so for instance $V_0(Q)=Q$, $V_1(Q)= Q \cup \{ x: x \text{ is a set with } \mathsf{rk}(x)=0\}$, and so on. 

We also stipulate to identify the urelements with the elements of $V(Q)$ of rank $<0$, for uniformity. 

\begin{Definition}\label{def:vdc}
    ($\mathsf{ZFU}$)
  We stipulate to extend the order $\leq_Q$ to a relation on $V(Q)$ as follows:
  \[
    x\leq_Q y \Leftrightarrow
    \begin{cases}
        x\leq_Q y & \text{ for } x,y\in Q\\
        x\in y & \text{ for } x\in Q,y\not\in Q\\
        x\subseteq y & \text{ for } x,y\not\in Q
    \end{cases}
  \]
  
  Let $\vdc Q$ be the class of all elements of $V(Q)$ that are either urelements from $Q$ or sets such that both them and all their set-elements are non-empty and downward-closed with respect to the extension of $\le_Q$ to $V(Q)$ we just described. 
  
  Formally,
  we say that $x\in V(Q)$ is \emph{downward-closed with respect to $\leq_Q$},
  and write $DC_{\leq_Q}(x)$, iff
  \[
    x\in Q \vee \forall y,z(y\in x \wedge z\leq_Q y \rightarrow z\in x).
  \]
  Then, we have that
  \[
    x\in \vdc Q \Leftrightarrow 
    DC_{\leq_Q}(x) \wedge 
    x\neq \emptyset \wedge
    \forall y\in x (DC_{\leq_Q}(y) \wedge y\neq\emptyset).
  \]

  For every ordinal $\alpha$, we put $\idcs \alpha Q$ to consist of all elements of $\vdc{Q}$ that are either urelements or sets of the rank $<\alpha$, i.e. $\idcs{\alpha}{Q}=\vdc{Q}\cap V_{\alpha}(Q)$. 
\end{Definition}

As we were saying, this Defintion is too complex for
our purposes: this is due to the universal quantification in the definition
of downward-closedness. This will be fixed in \ref{def:real-dcsets}.

We list some properties of $\vdc Q$ and $\idcs \alpha Q$, which are 
easily provable in, say, $\mathsf{ZFU}$.


\begin{Lemma}\label{lem:it_d-basic}
For any ordinal $\alpha$:
\begin{enumerate}
  \item sets in $\vdc{Q}$ ($\idcs{\alpha+1}{Q}$) are precisely the non-empty downward closed subsets of $\vdc{Q}$ ($\idcs{\alpha}{Q}$);
  \item all sets in $\vdc{Q}$ are transitive and their intersections with $Q$ are downward closed subsets of $Q$;
  \item $\vdc{Q}$ ($\idcs{\alpha}{Q}$) is a quasi-order;
  \item $\vdc{Q}$ ($\idcs{\alpha}{Q}$) is a partial order if $Q$ is a partial order;
  \item $\vdc{Q}$ ($\idcs{\alpha}{Q}$) is an end extension of any $\idcs{\beta}{Q}$ (any $\idcs{\beta}{Q}$, for $\beta<\alpha$).
 \end{enumerate}
\end{Lemma}

We now explore a different approach to define similar subclasses of $V(Q)$: 
in this case, we at first restrict $V(Q)$ in a very limited way, namely by
only taking care of the "non-emptiness" part. 

\begin{Definition}\label{def:real-dcsets}
  ($\prsou$) Let $\dot V(Q)$ consists of all elements of $V(Q)$ whose transitive closure does not contain the empty set. Let $\dot V_{\alpha}(Q)$ be the subclass of $\dot V(Q)$ consisting of elements of the rank $<\alpha$. 

  We define the class function $\supp: \dot V(Q) \to \mathcal{P}(Q)$ as
  \[
    \supp(x)=\{ x\} \text{ if } x\in Q \quad \text{ and } \quad
    \supp(x)= \bigcup \{\supp (x'): x'\in x\} \text{ otherwise}.
  \]
  $\supp(x)$ is called the \emph{support} of $x$.
\end{Definition}

  Let  $\dot{\mathcal{P}}$ be the operator mapping a set to the set of all its non-empty subsets. Then, an alternative definition of the hierarchy above (that would be suitable when we are working with set-sized orders in $\mathsf{ZFU}$) is 
  $$\dot V_{\alpha}(Q)=Q\cup\bigcup\limits_{\beta<\alpha}\dot{\mathcal{P}}(\dot V_{\beta}(Q)).$$

The trick is now to define a more convoluted way of comparing elements of these classes. 
As we will see, there are (at least) two sensible choices to do this.

\begin{Definition}\label{def:lesssim}
    ($\prsou$)
  The pre-orders $x \lesssim_Q^* y$ and $x \lesssim_Q y$ on $\dot V(Q)$ are defined by primitive-recursion:
  $$x\lesssim_Q^* y \iff \begin{cases} x\leq_Q y & \text{if }x,y\in Q\\ x\lesssim_Q^* \{y\} &\text{if }y\in Q \text{ and }x\not\in Q\\ \{x\}\lesssim_Q^*y  &\text{if }x\in Q \text{ and }y\not\in Q \\ \forall x'\in x\exists y'\in y(x'\lesssim_Q^*y') &\text{if }x,y\not\in Q\end{cases}$$ 
$$x\lesssim_Q y \iff \begin{cases} x\leq_Q y & \text{if }x,y\in Q\\ \text{always false} &\text{if }y\in Q \text{ and }x\not\in Q\\ \{x\}\lesssim_Q y  &\text{if }x\in Q \text{ and }y\not\in Q \\ \forall x'\in x\exists y'\in y(x'\lesssim_Qy') &\text{if }x,y\not\in Q.\end{cases}$$
In line with what we have done so far, we identify the sets 
$\dot V(Q)$, $\dot V_\alpha(Q)$, $\vdcalt{Q}$ and $\idcsalt{\alpha}{Q}$ with, respectively, the quasi-orders $(\dot V(Q),\lesssim_Q)$, $(\dot V_\alpha(Q), \lesssim_Q)$, $(\dot V({Q}),\lesssim^*_Q)$ and $(\dot V_{\alpha}({Q}),\lesssim^*_Q)$.
\end{Definition}

\begin{Remark}
  Recently, Anton Freund \cite{3-bqo-freund} used what in the terminology above is $\lesssim_Q^*$ on the set of hereditarily \emph{countable} sets $H_{\aleph_1}(Q)$ to show that the fact that the antichain with three elements is bqo implies $\aca^+$ over $\rca$. The restriction to hereditarily countable sets here is mostly because they are the kinds of sets that could be coded in the language of second-order arithmetic. A more important difference from $V^*(Q)$ is the lack of restriction on the presence of the empty set in the transitive closure. However, the impact of this on the order type is very controllable: It is easy to see that the difference between $(V_\alpha(Q),\lesssim_Q^*)/{\sim_Q^*}$ and $(\idcsalt{\alpha}{Q},\lesssim_Q^*)/{\sim_Q^*}$ is that the former has an additional well-ordered chain of the length $\alpha$ below all the equivalence class that have members from $\idcsalt{\alpha}{Q}$. For $(H_{\aleph_1}(Q),\lesssim_Q^*)/{\sim_Q^*}$ and $(H_{\aleph_1}(Q)\cap\vdcalt{Q},\lesssim_Q^*)/{\sim_Q^*}$ the added chain is of the length $\aleph_1$.
\end{Remark}

We got a bit ahead of ourselves in the Definition above: technically, we did
not verify that the relations $\lesssim_Q$ and $\lesssim^*_Q$ are transitive and 
reflexive, but it is immediate to do so by $\in$-induction.

Another obvious consideration is that if $x,y\in \dot V(Q)$ are such that 
$x\lesssim_Q y$ holds, then $x\lesssim^*_Q y$ holds as well.

Moreover, it is easy to see that all the considerations above would still hold if 
we were to restrict them to $\idcs {\alpha} Q$, $\dot V_{\alpha}(Q)$ and
$\idcsalt {\alpha} Q$, for every ordinal $\alpha$.

For the reader's convenience, we collect the facts we just stated in the following Lemma.

\begin{Lemma}
The following is provable in $\prsou$.
  \begin{enumerate}
  \item  $\lesssim_Q^*$ and $\lesssim_Q$ are qo's on $\dot V({Q})$;
  \item $\lesssim_Q^*$ is an extension of $\lesssim_Q$;
  \end{enumerate}
\end{Lemma}
\begin{proof}
    We show that $\lesssim_Q$ is a qo, to give an idea of how these sorts
    of proof work in $\prsou$.
    
    We need to see that $\lesssim_Q$ is reflexive and transitive. 
    We start with reflexivity, i.e.\ $\forall x ( x\in \dot V (Q)
    \rightarrow x\lesssim_Q x)$: we do this
    by straightforward $\in$-induction. Namely, if $x\in Q$, this follows 
    from reflexivity of $\leq_Q$. If $x\not\in Q$, we can
    suppose that $\forall y\in x(y\in \dot V(Q) \rightarrow y\lesssim_Q y)$.
    Since $x\in \dot V (Q)\wedge y\in x\rightarrow y\in \dot V(Q)$, it follows
    that $\forall x'\in x\exists y\in x(x'\lesssim_Q y)$ (namely, taking $x'=y$),
    so that $x\lesssim_Q y$.

    Now we show transitivity: let $x,y,z\in \dot V (Q)$ be such that 
    $x\lesssim_Q y \lesssim_Q z$, and let $a=\mathsf{TC}(\{x,y,z\})$. 
    We focus on the case that none of $x,y,z$ is a urelement: the other
    cases can be easily reconstructed from this one.
    We show by induction on the ordinals (hence, by $\in$-induction) that
    \[
    \forall \alpha (
    \forall x',y',z' \in a (
    \mathsf{rk}(x'),\mathsf{rk}(y'),\mathsf{rk}(z') \leq \alpha \rightarrow (
    x'\lesssim_Q y' \lesssim_Q z' \rightarrow x' \lesssim_Q z' ))).         
    \]
    Suppose the claim holds for every $\beta<\alpha$, and consider 
    $x',y',z'\in a$ such that $\mathsf{rk}(x'),\mathsf{rk}(y'),\mathsf{rk}(z') \leq \alpha$
    and $x'\lesssim_Q y' \lesssim_Q z'$. This means that 
    $\forall x''\in x'\exists y''\in y'(x''\lesssim_Q y'')$ and 
    $\forall y''\in y'\exists z''\in z'(y''\lesssim_Q z'')$: since the rank
    of these objects is less than $\alpha$, and since all of them are elements
    of $a$, we can conclude that $\forall x''\in x'\exists z''\in z'(x''\lesssim_Q z'')$,
    and hence that $x'\lesssim_Q z'$. This concludes the induction, and hence
    the proof, since $x,y,z\in a$.
\end{proof} 

We now explore the relationship between $\vdc Q$ and $\dot V(Q)$: 
it should be clear that $\vdc Q\subseteq \dot V(Q)$, 
and that in general $\dot V(Q)$ is strictly larger that $\vdc Q$. So it is 
natural to ask what happens if we take $x,y\in \vdc Q$ and compare them according to
$\lesssim_Q$: the following Lemma clarifies the situation. 

\begin{Lemma}\label{lem:leq-dc-notdc}
    ($\mathsf{ZFU}$) Suppose $Q$ is set-sized. Then, for any $x,y\in \vdc Q$, we have $x\leq_Q y$ iff $x\lesssim_Q y$.
    Moreover, for every $y\in \dot V (Q)$, there is an $x\in \vdc Q$ such that $y\sim_Q x$.
    (intuitively: the map $i:\vdc Q\to \dot V(Q)/\sim_Q$ associating $x\in \vdc Q$ to its equivalence class $[x]_{\sim_Q}$ is an order isomorphism) 
    
    The same holds for $\idcs \alpha Q$ and $\dot V_\alpha (Q)$.  
\end{Lemma}
\begin{proof}
    We only prove the result for $\vdc Q$ and $\dot V (Q)$, the relativization to sets of bounded rank being obvious.

    Suppose $x\leq_Q y$ holds. If $x,y\in Q$, we are done. If $x\in Q$ and 
    $y\not\in Q$,
    the assumption is that $x\in y$, and we want to show that 
    $\exists y'\in y (x\lesssim_Q y)$: but since we know that
    $x\in y$ and $x\lesssim_Q x$, we are done.    
    The case that $x,y\not\in Q$ and $x\subseteq y$ is similar:
    it is obvious that for every $x'\in x$ there is $y'\in y$
    with $x'\lesssim_Q y'$: namely, $x'$ itself. 

    Now suppose that $x,y\in \vdc Q$ and $x\lesssim_Q y$. We prove the claim
    by induction on the maximum of the rank of $x$ and $y$.
    If $x\in Q$ and 
    $y\in Q$, there is nothing to prove. If $x\in Q$ and $y\not \in Q$, 
    $x\lesssim_Q y$ means $\exists y'\in y(x\lesssim_Q y')$, so by induction
    on rank $x\in y'$, and by the fact that $y$ is a transitive set
    it follows that $x\in y$. Similarly, if $x\not\in Q$, 
    $\forall x'\in x\exists y'\in y(x'\lesssim_Q y')$ holds, so by
    inductive assumption $x'\leq_Q y'$, and again we conclude by 
    $\leq_Q$-downward-closedness of $y$.

    Next, we show that for every $y\in \dot V(Q)$, there is $x\in \vdc(Q)$ with $y\sim_Q x$, with the assumption that $Q$ is a set. We do this by induction on the rank of $y$. Suppose we have proved it for every ordinal less than $\alpha$, and assume that $y$ has rank $\alpha$. Hence, for every $y'\in y$, there is $x'\in \vdc(Q)$ such that $y'\sim_Q x$. Let $a$ be the set of these $x'$. Using the fact that $Q$ is a set, we can build the set $\{x''\in \idcs \alpha Q: \exists x'\in a(x''\leq_Q x')\}$, which is then the desired $x\in \vdc Q$.
\end{proof}

We now list several observations regarding the $\dot V_\alpha (Q)$ that can be
helpful to familiarize with them, and some other useful facts
that we will use in the rest of the paper.

\begin{itemize}
\item We note that it was essential to assume that $Q$ consisted of urelements, since they have special status within the constructions and otherwise we would be getting collisions with sets.
\item Although formally it would not be correct to talk about $\dot V_{\beta}({\dot V_{\alpha}({Q}))}$, since $\dot V_{\alpha}({Q})$ does not consists of urelements (for $\alpha>0$), one could still consider order on urelements $P$ isomorphic to $\idcs{\alpha}{Q}$ and consider then the order $\dot V_{\beta}({P})$ that morally should correspond to $\dot V_{\beta}({\dot V_{\alpha}({Q}))}$. One might have supposed that $\dot V_{\beta}({P})$ is isomorphic to $\dot V_{\alpha+\beta}({Q})$, however it is not in general the case. For example consider the case of $\alpha,\beta=1$ and  $Q$ being the one-element order $\{\star\}$
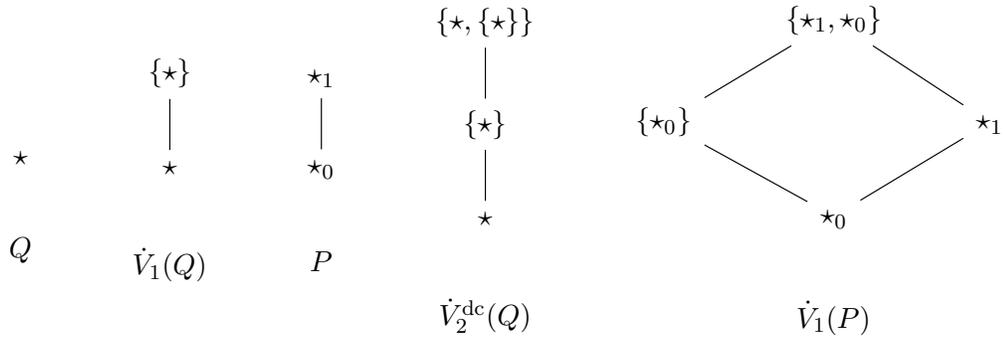
\begin{center}
  \begin{tikzcd}
   {} \arrow[d,phantom]\\
  \star \arrow[d,phantom]\\
  Q
\end{tikzcd}
\;\;\;\;\;\;\;
\begin{tikzcd}
  \{\star\}
  \arrow[d,dash] \\
  \star \arrow[d,phantom]\\
  \dot V_{1}({Q})
\end{tikzcd}
\;\;\;\;\;\;\;
\begin{tikzcd}
  \star_1
  \arrow[d,dash] \\ 
  \star_0 \arrow[d,phantom]\\
  P
\end{tikzcd}
\;\;\;\;\;\;\;
\begin{tikzcd}
  \{\star,\{\star\}\}
  \arrow[d,dash] \\ 
  \{\star\}
  \arrow[d,dash] \\
  \star \arrow[d,phantom]\\
  \idcs{2}{Q}
\end{tikzcd}
\;\;\;\;\;\;\;
\begin{tikzcd}
  &\{\star_1,\star_0\}\arrow[dl,dash]\arrow[dr,dash]&\\
  \{\star_0\}\arrow[dr,dash]&&\star_1\arrow[dl,dash]\\
  &\star_0\arrow[d,phantom]&\\
  &\dot V_{1}({P})&
\end{tikzcd}
\end{center}
\item To see that $\vdcalt Q$ behaves similarly consider the following example:

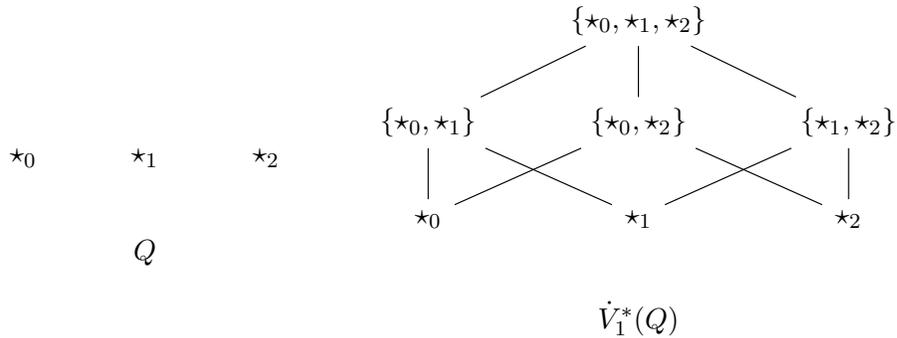
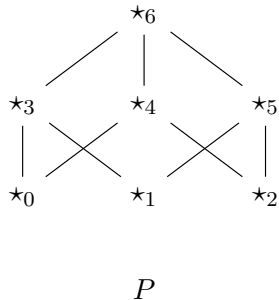
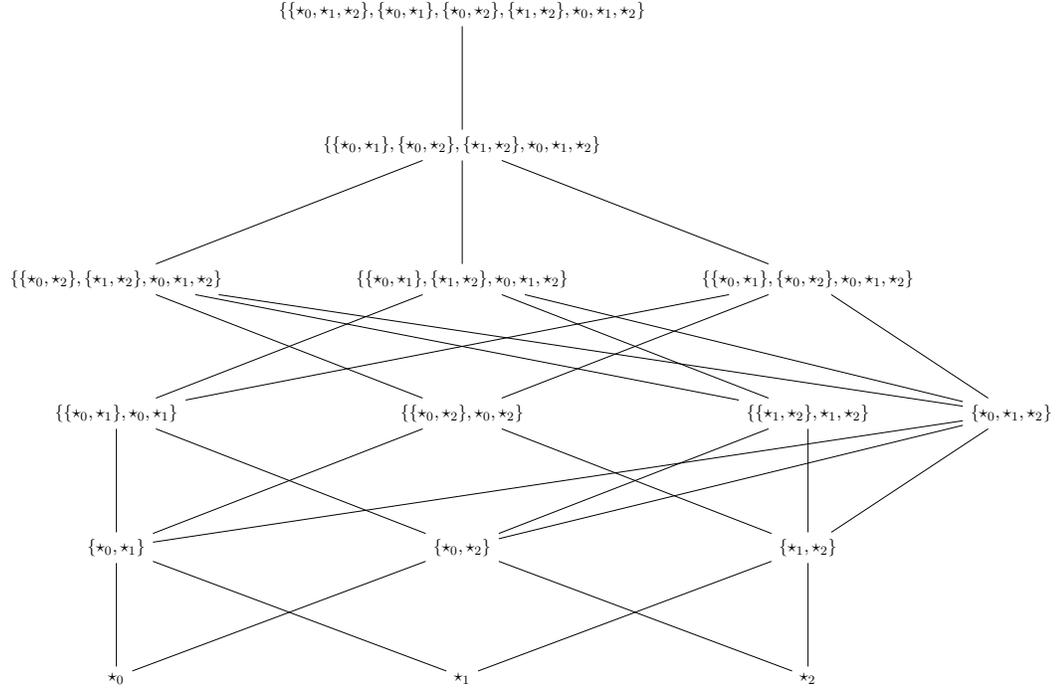
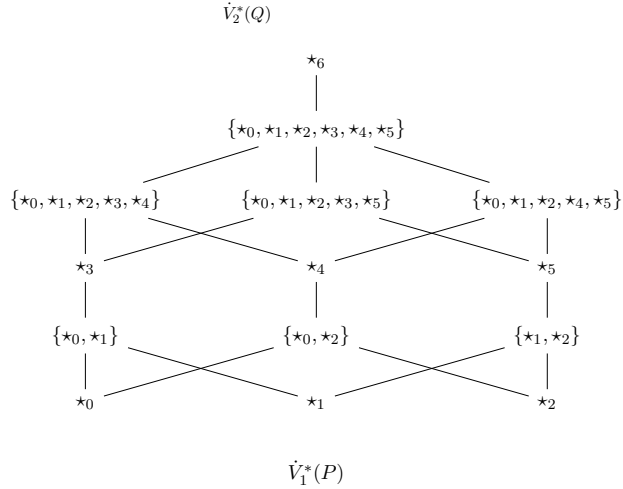
\begin{center}
\begin{tikzcd}
  & {} \arrow[d,phantom]& \\
  \star_0 \arrow[dr,phantom]& \star_1 &\star_2\\
  &Q&
\end{tikzcd}
\;\;\;\;\;\;\;
\begin{tikzcd}
  & \{\star_0,\star_1,\star_2\} \arrow[dl,dash]\arrow[d,dash]\arrow[dr,dash]& \\
   \{\star_0,\star_1\} \arrow[d,dash]\arrow[dr,dash]& \{\star_0,\star_2\} \arrow[dl,dash]\arrow[dr,dash] & \{\star_1,\star_2\} \arrow[dl,dash]\arrow[d,dash]\\
  \star_0 \arrow[dr,phantom]& \star_1&\star_2\\
  &\idcsalt{1}{Q}&
\end{tikzcd}
\end{center}
\;\;\;\;\;\;\;
\begin{center}
\begin{tikzcd}
  & \star_6 \arrow[dl,dash]\arrow[d,dash]\arrow[dr,dash]& \\
  \star_3\arrow[d,dash]\arrow[dr,dash]& \star_4\arrow[dl,dash]\arrow[dr,dash] &\star_5\arrow[dl,dash]\arrow[d,dash] \\ 
  \star_0 \arrow[dr,phantom]& \star_1 &\star_2\\
  &P&
\end{tikzcd}\\
\end{center}

\begin{center}
\adjustbox{scale=0.6,center}{
\begin{tikzcd}[row sep=6em]
  &\{\{\star_0,\star_1,\star_2\},\{\star_0,\star_1\},\{\star_0,\star_2\},\{\star_1,\star_2\},\star_0,\star_1,\star_2\} \arrow[d,dash] & & \\
  &\{\{\star_0,\star_1\},\{\star_0,\star_2\},\{\star_1,\star_2\},\star_0,\star_1,\star_2\} \arrow[dl,dash]\arrow[d,dash]\arrow[dr,dash] & & \\
  \{\{\star_0,\star_2\},\{\star_1,\star_2\},\star_0,\star_1,\star_2\} \arrow[dr,dash]\arrow[drr,dash]\arrow[drrr,dash]& \{\{\star_0,\star_1\},\{\star_1,\star_2\},\star_0,\star_1,\star_2\} \arrow[dl,dash]\arrow[dr,dash]\arrow[drr,dash]& \{\{\star_0,\star_1\},\{\star_0,\star_2\},\star_0,\star_1,\star_2\} \arrow[dll,dash]\arrow[dl,dash]\arrow[dr,dash] & \\
   \{\{\star_0,\star_1\},\star_0,\star_1\}  \arrow[d,dash]\arrow[dr,dash]& \{\{\star_0,\star_2\},\star_0,\star_2\} \arrow[dl,dash]\arrow[dr,dash] & \{\{\star_1,\star_2\},\star_1,\star_2\}  \arrow[dl,dash]\arrow[d,dash] & \{\star_0,\star_1,\star_2\}  \arrow[dlll,dash]\arrow[dll,dash]\arrow[dl,dash]\\
   \{\star_0,\star_1\} \arrow[d,dash]\arrow[dr,dash]& \{\star_0,\star_2\} \arrow[dl,dash]\arrow[dr,dash] & \{\star_1,\star_2\} \arrow[dl,dash]\arrow[d,dash] &\\
  \star_0 \arrow[dr,phantom]& \star_1&\star_2 &\\
  &\idcsalt{2}{Q}& &
\end{tikzcd}
}
\end{center}

\begin{center}
\adjustbox{scale=0.7}{
\begin{tikzcd}
    &\star_6 \arrow[d,dash]& \\
  & \{\star_0,\star_1,\star_2,\star_3,\star_4,\star_5\}  \arrow[dl,dash]\arrow[d,dash]\arrow[dr,dash] &\\
   \{\star_0,\star_1,\star_2,\star_3,\star_4\}  \arrow[d,dash]\arrow[dr,dash]& \{\star_0,\star_1,\star_2,\star_3,\star_5\}  \arrow[dl,dash]\arrow[dr,dash]& \{\star_0,\star_1,\star_2,\star_4,\star_5\}  \arrow[dl,dash]\arrow[d,dash]\\
   \star_3\arrow[d,dash]& \star_4\arrow[d,dash] &\star_5\arrow[d,dash] \\
    \{\star_0,\star_1\} \arrow[d,dash]\arrow[dr,dash]& \{\star_0,\star_2\} \arrow[dl,dash]\arrow[dr,dash] & \{\star_1,\star_2\} \arrow[dl,dash]\arrow[d,dash] \\
  \star_0 \arrow[dr,phantom]& \star_1 &\star_2\\
  &\idcsalt{1}{P}&
\end{tikzcd}
}
\end{center}

Here we omitted equivalent elements, e.g. in $\idcsalt{1}{Q}$, we have $\star_0\mathrel\sim\{\star_0\}$.
\end{itemize}

\subsection{The sets $\dot V^{dc,f}_\alpha (Q)$} 
We will introduce one last approach to the iterated non-empty downward-closed sets, obtaining the classes and $\dot V^{dc,f}(Q)$ and $\dot V^{dc,f}_\alpha (Q)$. But to do this, we make one extra assumption, that will hold for the present Subsection only: that the class of urelements $Q$ is actually a set.

\begin{Definition}
    ($\prsou$ $+$ ``$Q$ is a set")
    Let $\leq_Q$ be the order on $V(Q)$ defined in \Cref{def:vdc}. We define
     the primitive recursive function $D^f_Q \colon V(Q) \to V(Q)$  as
     \[
        D^f_Q(P)= \{ \{p\in P : \forall a\in A (a \not\leq_Q p )\} \colon
        A \in [P]^{<\omega}\}.
     \]
     We then define the sets $\dot V^{dc,f}_\alpha (Q)$ by primitive
     recursion by putting 
     \[
        \dot V^{dc,f}_\alpha (Q) = Q \cup 
        \bigcup_{\beta<\alpha} D^f_Q(\dot V^{dc,f}_\beta (Q))
        \cap \dot V(Q).
     \]
     We let $\dot V^{dc,f} (Q)= \bigcup_{\alpha}\dot V^{dc,f}_\alpha (Q)$.

     Again, we identify $\dot V^{dc,f}_\alpha (Q)$ and $\dot V^{dc,f} (Q)$ with the orders
     $(\dot V^{dc,f}_\alpha (Q),\leq_Q)$ and $(\dot V^{dc,f} (Q),\leq_Q)$,
     respectively.
\end{Definition}

So, the fact that we restrict to the \emph{finite} powerset gives us a 
hierarchy of sets. We highlight that the element of $\dot V^{dc,f}_\alpha (Q)$ are 
downward-closed in $\dot V^{dc,f}_\alpha (Q)$ with respect to the order $\leq_Q$: if $x,y,z\in \dot V^{dc,f}_\alpha (Q)$ and $y\in x$ and $z\leq_Q y$ both hold, then $z\in x$ holds as well, as one easily derives from the definition: intuitively, this follows from the fact that the elements of $\dot V^{dc,f}_\alpha (Q)$ are complement of upward-closed sets.

It is apparent from the definition that every $\dot V^{dc,f}_\alpha (Q)$ is a subset
of $\dot V_{\alpha}(Q)$. As in the case of $\vdc Q$, the order is respected in the restriction.

\begin{Lemma}
    ($\prsou$$+$ ``$Q$ is a set")
    Let $x,y\in \dot V^{dc,f} (Q)$. Then, $x\leq_Q y$ iff $x\lesssim_Q y$.
\end{Lemma}
\begin{proof}
    The proof is obviously very similar to that of \ref{lem:leq-dc-notdc}. In particular,
    the proof that $x\leq_Q y$ implies $x\lesssim_Q y$ can be copy-pasted from there.

    Let us now suppose that $x\lesssim_Q y$. We prove by induction that 
    \[
        \forall \alpha (\forall x,y\in \dot V^{dc,f}_\alpha (Q)(
        x\lesssim_Q y \rightarrow x\leq_Q y)).
    \]
    Suppose the claim holds for all $\beta<\alpha$, and we prove it for $\alpha$. 
    If $x,y\in Q$, there is nothing to prove. If $x\in Q$ and $y\not\in Q$, then
    $x\lesssim_Q y$ means that $x\lesssim_Q y'$ for $y'\in y$, and we conclude by inductive assumption. Finally, if $x,y\not\in Q$, $x\lesssim_Q y$ means that $\forall x'\in x \exists y'\in y(x'\lesssim_Q y')$, which by inductive assumption implies $\forall x'\in x \exists y'\in y(x'\leq_Q y')$: since we have already observed that elements of $\dot V^{dc,f}_\alpha (Q)$ are downward-closed in $\dot V^{dc,f}_\alpha (Q)$, we conclude that $x\subseteq y$, and so $x\leq_Q y$. 
\end{proof}

There is a close relationship between $\dot V^{dc,f}_\alpha (Q)$
and $\dot V_\alpha (Q)$, if any of the two is well-founded.

\begin{Lemma}\label{lem:vdc-wf-skel}
    ($\prsou$ $+$ ``$Q$ is a set")
    Suppose that $\dot V^{dc,f}_\alpha(Q)$ is well-founded. Then, for every $y\in \dot V_\alpha (Q)$, there is an $x\in \dot V^{dc,f}_\alpha (Q)$ such that $y\sim_Q x$
    (intuitively, the map $i:\dot V^{dc,f}_\alpha Q\to \dot V(Q)/\sim_Q$ associating $x\in  Q$ to its equivalence class $[x]_{\sim_Q}$ is an order isomorphism). 
\end{Lemma}
\begin{proof}
    First, we show that if $\dot V^{dc,f}_\alpha(Q)$ is well-founded, then for every $\beta<\alpha$ $\dot V^{dc,f}_\beta(Q)$ is wqo. We prove this by cases: if $\alpha=0$, the claim is trivial. If $\alpha=\beta+1$, we proceed by contradiction and suppose that there is a bad sequence $(x_i)_{i\in\omega}$ in $\dot V^{dc,f}_\beta(Q)$: then we can build the sequence $(y_i)_{i\in\omega\setminus 0}$ of elements of $\dot V^{dc,f}_\beta(Q)$ defined as
    \[
        y_i=\{ x\in \dot V^{dc,f}_\beta(Q): \forall j< i(x\not\leq_Q x_j)\},
    \]
    and notice that it is descending. Finally, if $\alpha$ is limit, $\dot V^{dc,f}_\alpha(Q)$ is just the union of the previous levels, so if any $\dot V^{dc,f}_\beta(Q)$ is not wqo, $\dot V^{dc,f}_{\beta+1}(Q)$ is not well-founded, and so neither is $\dot V^{dc,f}_\alpha(Q)$.

    Next, we prove by $\in$-induction 
    \[
        \forall y(y\in \dot V_\alpha(Q)\rightarrow \exists x\in \dot V^{dc,f}_\alpha(Q)
        (x\sim_Q y)).
    \]
    If $y\in Q$, there is nothing to show, so suppose $y$ is a set. We define
    \[
        x=\{ x'\in \dot V^{dc,f}_\alpha(Q): \exists y'\in y (y'\sim_Q x')\}.
    \]
    It is immediate that $x\lesssim_Q y$. Moreover, by inductive assumption, for every $y'\in y$ there is $x'\in \dot V^{dc,f}_\alpha(Q)$ with $x'\sim_Q y'$, so $y\lesssim_Q x$ holds as well. Finlly, to prove that $x\in \dot V^{dc,f}_\alpha(Q)$, notice that $x$ is a downward-closed subset of $\dot V^{dc,f}_\beta(Q)$ for some $\beta<\alpha$, and hence it is finitely generated by the well-known properties of wqo's.
\end{proof}

The fact that $\dot V^{dc,f}_\alpha(Q)$ is or is not well-founded depends on the degree of wqo-ness of the set $Q$ of urelements, as we will see in the next Subsection.

\subsection{Unwinding of sequences and winding of arrays.}\label{subsec:wind-unwind}

In this subsection, we see how to move between bad sequences into $\dot V_\alpha(Q)$ and bad arrays of rank $\leq\alpha$ into $Q$.

\begin{Lemma}\label{lem:unwinding}
($\prsou$)
If $f\colon \omega\to \dot V_\alpha (Q)$ (resp. $f^*\colon \omega\to \idcsalt\alpha Q$) is a bad sequence for $\lesssim_Q$  ($\lesssim_Q^*$), then there is a bad array $g\colon F\to Q\times \omega$ ($g^*\colon F^*\to Q$) of the rank $\le \alpha$. 
\end{Lemma}
Before giving the proof, we highlight two things worth remarking: the first
is that in the proof, it will become apparent why we were careful in keeping
the empty set out of $\dot V_{\alpha}(Q)$ in the definition above; the second is
that the elegant form of the result is the main reason why went for an
``unorthodox'' definition of $\alpha$-wqo-ness, instead of sticking with the
classical one in terms of lexicographic rank of fronts.

\begin{proof}
  Notice that the claim is obvious for $\alpha=0$, so can suppose that
  $\alpha$ is non-zero.

  We will prove the claim by contradiction: we suppose that there is a bad
  sequence $f:\omega\to \dot V_\alpha (Q)$ (resp. $f^*\colon \omega\to \idcsalt\alpha Q$), and ``unravel'' it to produce a bad
  array $g:F\to Q\times\omega$, ($g^*:F^*\to Q$) where $F$ ($F^*$) is a front with $\rank F\leq \alpha$.

 In both the case of $\dot V_\alpha (Q)$ and $\idcsalt \alpha Q$, the initial step is the construction of auxiliary functions $h$ and $h^*$ whose domain is $[\omega]^{<\omega}\setminus \{\emptyset\}$. We start with the case of $\idcsalt \alpha Q$, which is simpler and gives a clearer intuition of what is going on.

  We define values $h^*(n_0,\ldots,n_{k-1})$ of the function $h^*:[\omega]^{<\omega}\setminus \{\emptyset\} \to \idcsalt \alpha Q$ by recursion on the length $k$ of the sequence $(n_0,\ldots,n_{k-1})\in [\omega]^{<\omega}$. In the following clauses we go to the next clause if all the previous ones are non-applicable:
  \begin{enumerate}
  \item $h^*(n)=f(n)$ if $f(n)\in Q$;
  \item $h^*(n)=f(n)$ if $f(n)$ is a set;
  \item \label{hs_case_3} $h^*(n_0,\dots,n_{k})=q$ if $h^*(n_0,\dots,n_{k-1})=q\in Q$;
  \item \label{hs_case_4} we set $h^*(n_0,\dots,n_k)=q\in Q$, for a $q\in Q$ such that $q\in \supp( h^*(n_0,\dots,n_{k-1}))$ and $q\not\lesssim^*_Q h^*(n_1,\dots,n_{k})$, if such a $q$ exists.
  \item \label{hs_case_5} we set $h^*(n_0,\dots,n_k)=x$ for some $x\in h^*(n_0,\dots,n_{k-1})$ such that for every $y\in h^*(n_1,\dots,n_k)$, $x\not\lesssim^*_Q y$ holds.
 \end{enumerate}
 Note that in order for this definition to be sound (i.e. that we would be able to pick $x$ in the last item) it is sufficient to guarantee that $h^*(n_0,\ldots,n_{k-1})\not \lesssim^*_Q h^*(n_1,\ldots,n_{k})$, i.e. that $h^*\upharpoonright_{[\omega]^k}$ is a bad array.

 We prove this by induction on $k$. If $k=1$, the claim is obvious. Suppose the claim is proved for a certain $k$, and we prove it for $k+1$. We show that $h^*(n_0,\dots,n_{k})\not \lesssim^*+Q h^*(n_1,\dots,n_{k+1})$ by cases. 
 \begin{itemize}
     \item Suppose $h^*(n_0,\dots,n_{k})$ was defined by (\ref{hs_case_3}). Since $h^*(n_0,\dots,n_{k-1})\not\lesssim^*_Q h^*(n_1,\dots,n_{k})$ and $h^*(n_1,\dots,n_{k+1})\lesssim^*_Qh^*(n_1,\dots,n_{k})$, we are done.
     \item Suppose $h^*(n_0,\dots,n_{k})$ was defined by (\ref{hs_case_4}). We are done by noting that $h^*(n_1,\dots,n_{k+1})\lesssim^*_Qh^*(n_1,\dots,n_{k})$ and that $q\in \supp(h^*(n_0,\dots,n_{k-1}))$ implies $q\lesssim^*_Q h^*(n_0,\dots,n_{k-1})$.
     \item Finally, suppose $h^*(n_0,\dots,n_{k})$ was defined by (\ref{hs_case_5}). Suppose at first that $h(n_1,\dots,n_k)\in Q$: by \cite[Remark 2.3]{3-bqo-freund}, there has to be $q\in Q$ such that $q\in \supp( h^*(n_0,\dots,n_{k-1}))$ and $q\not\lesssim^*_Q h^*(n_1,\dots,n_{k})$. 
     More explicitly, we can show by induction that $\forall x\forall p\in Q(\forall q\in \supp(x) (q \lesssim^*_Q p)\rightarrow x\lesssim^*_Q p)$. Indeed, the claim is true if $x\in Q$. Suppose it is not, then by induction we have that $\forall x'\in x(x\lesssim^*_Q p)$, which implies $x\lesssim^*_Q \{p\}\sim^*_Q p$, as we wanted.
     So we should have used (\ref{hs_case_4}), and so this case cannot occur. 
     
     So suppose that $h^*(n_1,\dots,n_{k})$ is a set: then $h^*(n_1,\dots,n_{k+1})\in h^*(n_1,\dots,n_{k})$, and so by definition $h^*(n_0,\dots,n_{k})\not\lesssim^*_Q h^*(n_1,\dots,n_{k+1})$
 \end{itemize}

 We now move to the definition of $h$ in the case of $\dot V_\alpha (Q)$.

  Here, we will employ a slightly larger order $Q_\alpha$ with domains $(\dot V_\alpha Q \setminus Q)\cup (Q\times \omega)$, where we compare $x,y\in \dot V_\alpha (Q) \setminus Q$ as before, $x,y\in Q\times \omega$ as the elements of product of the product order, and $(x,n)\in Q\times \omega$ with $y\in \dot V_\alpha (Q) \setminus Q$ as we would compare $x$ with $y$ in $\dot V_\alpha (Q)$. We denote this order on $Q_\alpha$ by $\leq_\alpha$.

  
  We define values $h(n_0,\ldots,n_{k-1})$ of the function $h$ by recursion on the length $k$ of the sequence $(n_0,\ldots,n_{k-1})\in [\omega]^{<\omega}$. In the following clauses we go to the next clause if all the previous ones are non-applicable:
  \begin{enumerate}
  \item $h(n)=(f(n),0)$ if $f(n)\in Q$;
  \item $h(n)=f(n)$ if $f(n)$ is a set;
  \item \label{h_case_3} $h(n_0,\ldots,n_k)=h(n_0,\ldots,n_{k-1})$ if $h(n_0,\ldots,n_{k-1})\in Q\times \omega$;
  \item \label{h_case_4} $h(n_0,\ldots,n_k)=(q,k)$, where $q$ is picked to be any urelement
  such that $q\in \supp( h(n_0,\ldots,n_{k-1}))$ and $q\not\lesssim_Q h(n_1,\ldots,n_{k})$, if such a $q$ exists.
  Note that this is in particular always applicable if $h(n_1,\ldots,n_{k})\in Q\times \omega$, since the relation $\lesssim_Q$ is not defined on those elements;
  \item \label{h_case_5} we set $h(x_0,\dots,x_{k})=x$, for some set $x\in h(n_0,\ldots,n_{k-1})$ such that, for every $y\in h(n_1,\ldots,n_{k})$, $x\not\lesssim_Q y$.
  \end{enumerate}
Also in this case, to show that $h$ defines a total function, i.e.\ that at least one of the cases above is triggered, we only need to show that $h\rst_{[\omega]^{k}}$ is a bad array with respect to $\leq_\alpha$.  

  We prove the latter by induction on $k$. For $k=1$ we have this immediately from the badness of $f$. We will justify the step of induction by proving that $h(n_0,\ldots,n_{k})\not \leq_\alpha h(n_1,\ldots,n_{k+1})$ under the assumption that $h(n_0,\ldots,n_{k-1})\not \leq_\alpha h(n_1,\ldots,n_{k})$. To prove this we consider cases for the defining clause of of $h(n_0,\ldots,n_k)$:
\begin{itemize}
\item If $h(n_0,\ldots,n_k)$ was defined by (\ref{h_case_3}), then $h(n_0,\ldots,n_{k-1})=(q,l)\in Q\times \omega$. Since $(q,l) \not \leq_\alpha h(n_1,\ldots,n_k)$ by inductive assumption, and since it is always the case that $h(n_1,\ldots,n_k)\geq_\alpha h(n_1,\ldots,n_{k+1})$ (this follows easily from the definition), we see that $h(n_0,\ldots,n_{k})=(q,l)\not\leq_\alpha h(n_1,\ldots,n_{k+1})$. 
\item Assume $h(n_0,\ldots,n_k)=(q,k)\in Q\times \omega$ was defined by (\ref{h_case_4}). If $h(n_1,\ldots,n_k)$ is of the form $(p,l)\in Q\times \omega$, then $l<k$ and thus $(q,k)\not\leq_\alpha (p,l)$. If $h(n_1,\ldots,n_k)$ was a set, then $q\not\lesssim_Q y$ for any $y\lesssim_Q h(n_1,\ldots,n_k)$, and thus $(q,k)\not\leq_\alpha h(n_1,\ldots,n_{k+1})$.
\item Finally assume that $h(n_0,\ldots,n_k)=x\in \dot V_\alpha (Q) \setminus Q$ was defined by (\ref{h_case_5}). The value $h(n_1,\ldots,n_{k+1})$ is either $y\in h(n_1,\ldots,n_{k})$ or $(p,k)$ for $p\in \supp(h(n_1,\ldots,n_{k}))$. In the first case $x\not \leq_\alpha y$ by definition of $h(x_0,\dots,x_k)$. In the second case $x\not \le (p,k)=h(n_1,\ldots,n_{k+1})$ simply by definition of $Q_\alpha$.
\end{itemize}

Finally, we move to define the bad arrays $g$ and $g^*$, using the functions $h$ and $h^*$ we have just defined. We only do this in the case of $g$, the other one being virtually identical.
  
  Let us define the set $H:=\{ \sigma \in [\omega]^{<\omega}:
  h(\sigma)\in Q\times \omega\}$. Observe that
  the set $F$ is a front of the rank $\le \alpha$:  $$F:=\min_\sqsubseteq H=\{ \sigma\in H: \forall \sigma'\sqsubset
  \sigma (\sigma'\not\in H)\}.$$
  
  Indeed, this is witnessed by the ranking function $e:T_F\to \alpha +1$ such that $e(\emptyset)=\alpha$ and $e(n_0,\ldots,n_{k-1})=\mathsf{rk}(h(n_0,\ldots,n_{k-1}))$ (i.e., $(F,e)$ is a front).
  
  Clearly, $g=h\upharpoonright_F$ is an array with values from $Q\times \omega$. To show that it is bad, for any $\sigma,\tau\in F$, $\sigma\triangleleft \tau$, we fix $\sigma'=\sigma\cup \tau$ and $\tau'=((\sigma\cup \tau)\setminus \{\min(\sigma)\})\cup \{\max(\tau\cup \sigma)+1)\}$. Observe that $h(\tau)=h(\tau')$, $h(\sigma)=h(\sigma')$, and that since the lengths of $\sigma'$ and $\tau'$ are the same, we already proved that $h(\sigma')\not\le h(\tau')$.
\end{proof}

\begin{Remark}\label{remark:supp}
    Notice that, in the construction above, the first component of $g^*(\sigma)$ belongs to $\supp (f^*(\sigma(0)))$, and the same applies to $g$ and $f$, by examining the construction of $h^*$ and $h$.
\end{Remark}

Notice that the Lemma has the following corollary: if $Q$ is $(\alpha+1)$-wqo, then $\dot V_\alpha (Q)$ and $\idcsalt \alpha Q$ are wqo's, provably in $\prsou$: this is immediate in the case of $\idcsalt \alpha Q$, whereas in the other case we just have to notice that $\prsou$ proves that well-orders (and so in particular $\omega$) are bqo's, and that the product of $\alpha$-wqo's is $\alpha$-wqo.

\begin{Proposition}
($\prsou$)
If $g\colon F\to Q$  is a bad array  of the rank $\le \alpha$, then there is a bad sequence $f\colon \bigcup F\to  \idcsalt \alpha Q$.
\end{Proposition}
\begin{proof}
  The main idea of the proof is to do the opposite of what we did in
  the previous Lemma: whereas in that case we were ``unwinding''
  a sequence $f$ to find a certain array $g$,
  here we will go the other way, collecting elements
  in the range of a given array $g:F\to Q$ until we reach a sequence $f$
  with range $\idcsalt \alpha Q$, in a way "winding" the array. 

  We define an auxiliary function $h:[\bigcup F]^{<\omega} \to \dot V(Q)$ by recursion on the rank of $T_F$ as follows:
  \[
    h(\sigma)=
    \begin{cases}
        g(\tau) & \text{ if } \tau\sqsubseteq \sigma \text{ and } \tau\in F \\
        \{ h(\sigma \conc n) : \sigma\conc n \in F \} & \text{ otherwise.}
    \end{cases}
  \]
  The definition makes sense since, for every $\sigma\in [\bigcup F]^{<\omega}$, if $\sigma\not\in T_F$, then there is $\tau\in F$ with $\tau\sqsubseteq\sigma$.

  We start showing that the range of $h$ is actually $\dot V_{\alpha+1}(Q)$. To do this, we show that $\forall \sigma\in T_F (h(\sigma)\in \dot V_{\mathsf{rk}_F(\sigma)+1}(Q))$. We do this by induction on the rank of $\sigma$, namely we prove 
  \[
    \forall \beta (\forall \sigma \in T_F (\mathsf{rk}_F(\sigma)\leq \beta \rightarrow h(\sigma)\in \dot V_{\mathsf{rk}_F(\sigma)+1}(Q))
  \]
  If $\mathsf{rk}_F(\sigma)=0$, $h(\sigma)$ is a subset of $Q$, so the claim holds. Suppose that we have proved the claim for every $\gamma<\beta$, and let $\sigma\in T_F$ have rank $\beta$: since $h(\sigma\conc n)\in \dot V_{\mathsf{rk}_F(\sigma\conc n)+1}(Q)$ for every $n\in \bigcup F$ and $\mathsf{rk}_F(\sigma\conc n)+1\leq \mathsf{rk}_F(\sigma)$, it follows that $h(\sigma)\subseteq \dot V_{\mathsf{rk}_F(\sigma)}$, and so $h(\sigma)\in \dot V_{\mathsf{rk}_F(\sigma)+1}(Q)$.

  We let $f$ be $h\rst_{[\omega]^1}: \bigcup F\to \dot V_\alpha(Q)$.  We claim that $f$ is a bad sequence in $\idcsalt \alpha Q$ if $g$ is a bad array in $Q$.

  To do so, we first show that for every $n\geq 1$, if $\sigma,\tau\in [\bigcup F]^n$ are such that $\sigma\barprec \tau$ and $h(\sigma)\lesssim^*_Q h(\tau)$, then there is $m\in \bigcup F$ such that $h(\sigma\cup\tau)\lesssim^*_Q h(\tau\conc m)$.

  We show this by cases. If $h(\sigma),h(\tau)\in Q$, we can pick as $m$ any element of $\bigcup F$ larger than $\max \tau$. Similarly, if $h(\sigma)\not\in Q$ but $h(\tau)\in Q$, $h(\sigma)\lesssim^*_Q h(\tau)$ means that $\forall x'\in h(\sigma)(x'\lesssim^*_Q h(\tau))$, so again $m$ can be any element of $\bigcup F$ larger than $\max \tau$. If $h(\sigma)\in Q$ and $h(\tau)\not\in Q$, $h(\sigma)\lesssim^*_Q h(\tau)$ means that $\exists y'\in h(\tau) (h(\sigma)\lesssim^*_Q y')$: by definition of $h$, there is $m\in \bigcup F$ with $h(\tau\conc m)= y'$, and this is the $m$ we are after. Similarly, if $h(\sigma),h(\tau)\not\in Q$, $h(\sigma)\lesssim^*_Q h(\tau)$ means that $\forall x'\in h(\sigma)\exists y'\in h(\tau)(x'\lesssim^*_Q y')$, and again we pick the $m$ such that $h(\tau\conc m)=y'$.

  Suppose for a contradiction that $f$ was good: then there are $n_0<n_1\in \bigcup F$ such that $h(n_0)\lesssim^*_Q h(n_1)$. The claim we just proved allows us to build by recursion a set $X\subseteq [\bigcup F]^{\omega}$ such that for every $n$ $h(X\rst_n)\lesssim^*_Q h(X^-\rst_{n})$. Let $\sigma\in F$ be an initial segment of $X$, and $\tau\in F$ be an initial segment of $X^-$: then $g(\sigma)\leq_Q g(\tau)$, which gives the desired contradiction.
\end{proof}

We could rephrase this result as follows: if $\idcsalt \alpha Q$ is wqo, then $Q$ is $(\alpha+1)$-wqo. Since $\idcsalt \alpha Q$ is an extension of $\dot V_\alpha(Q)$, this implies that if $\dot V_\alpha(Q)$ is wqo, then $Q$ is $(\alpha+1)$-wqo.

\section{Connections with the transfinite sequences}\label{sec:tr-seq}

In this Section, we see that the classes $\dot V_\alpha (Q)$ (resp.\ $\idcsalt \alpha Q$) are strongly related to the qo of transfinite sequences up to a certain length ordered by (weak) embeddability. To do this, we will define maps that allow to move from one setting to the other. The indecomposable sequences will be an important ingredient for these maps, effectively serving as a skeleton for the structure of the larger order that the transfinite sequences form.

\begin{Definition}
($\prsou(\mathsf{enu})$)
Let $\iota$  be a  class function mapping each $\dot V_\alpha (Q)$ to $\seq {\omega^{1+\alpha}} Q$:
  \begin{enumerate}
  \item $\iota(q)=q$, for $q\in Q$.
  \item $\iota(x)=\sum_{i<\omega}  \iota(\mathsf{enu}(x,h(i)))$, for $x$ that are not urelements.
  \end{enumerate}
  Here $h\colon \omega\to \omega$ is some fixed function such that the $h$-preimage of any natural is infinite. For definiteness we put $h(x)$ to be the remainder of $x$ under the division by $[\sqrt{x+1}]$.
\end{Definition}

It is clear from the definition that $\iota (x)$ is indecomposable, and hence also weakly indecomposable, for every $x\in \dot V_\alpha(Q)$, due to the infinite repetition of the shorter subsegments of $\iota (x)$.

\begin{Lemma}
  For any $\alpha$, the map $\iota$ restricted to $\dot V_\alpha (Q)$ ($\idcsalt \alpha Q$) is an order-preserving and order-reflecting map with range $\indseq {\omega^{1+\alpha}} Q$ ($\indseqalt {\omega^{1+\alpha}} Q$).
\end{Lemma}
\begin{proof}
    We start by proving that the range is correct, namely that if $x\in \dot V_\alpha(Q)$, then $\iota(x)\in \seq {\omega^{1+\alpha}} Q$. For any $x\in \dot V(Q)$, let $a$ be a set such that $x\in a$ and $a=\mathsf{TC}(a)$. We prove by $\in$-induction that
    \[
        \forall \alpha (\forall x'\in a (x'\in \dot V_\alpha (Q) \rightarrow \iota(x)\in \seq {\omega^{1+\alpha}} Q)).
    \]
    Suppose the claim holds for $\beta<\alpha$ and we prove it for $\alpha$. If $\alpha$ is a successor or $0$, the claim is immediate. If instead $\alpha$ is limit, notice that $x'\in \dot V_\alpha(Q)$ implies $x'\in \dot V_\beta(Q)$ for $\beta<\alpha$, and the claim follows.

    Next, we show that $\iota$ is order-preserving. Given $x,y\in \dot V(Q)$, let $a$ be such that $x,y\in a$ and $a=\mathsf{TC}(a)$. We show simultaneously by induction that
    \[
        \forall \alpha (\forall x',y'\in a (x',y'\in \dot V_\alpha (Q) \wedge x'\lesssim_Q y' \rightarrow \iota(x') \preceq \iota(y'))) 
    \]
    and
    \[
        \forall \alpha (\forall x',y'\in a (x',y'\in \dot V_\alpha (Q) \wedge x'\lesssim^*_Q y' \rightarrow \iota(x') \preceq^* \iota(y'))). 
    \]
    Suppose the claims hold for all $\beta<\alpha$, and we show it for $\alpha$. We proceed by cases.
    \begin{itemize}
        \item If $x',y'\in Q$, both claims are obvious, since both $\iota(x')$ and $\iota(y')$ are sequences of length $1$ and equal to $x'$ and $y'$.
        \item If $x'\in Q$ and $y'\not\in Q$, there is $y''\in y'$ such that $x'\lesssim_Q y''$ (resp.\ $x'\lesssim^*_Q y''$), and the (weak) embedding is given by the inductive assumption.
        \item If $x',y'\not\in Q$, then for every $x''\in x'$ there is $y''\in y'$ with $x''\lesssim_Q y''$ (resp.\ $x''\lesssim^*_Q y''$. By inductive assumption, there is a (weak) embedding between $\iota(x'')$ and $\iota(y'')$ for every such $x''$ and $y''$. One can then combine these (weak) embeddings into a (weak) embedding of $\iota(x')$ into $\iota(y')$, using the fact that every $\iota(y'')$ appears as a subsegment in $iota(y')$ infinitely many times.
    \end{itemize}
    This concludes the proof for the non-weak case. In the weak case, we also have 
    \begin{itemize}
        \item If $x'\not\in Q$ and $y'\in Q$, then for all $x''\in x'$ $x''\lesssim^*_Q y'$, so by assumption there are weak embeddings of each $\iota(x'')$ into $\iota(y')$, all witnessed by the map sending every element of the domain to $0$. Such a map then also witnesses that $\iota(x')\preceq^* \iota(y')$.
    \end{itemize}
    The argument above entails that for the chosen $x,y$, if $x\lesssim_Q y$ (resp.\ $x\lesssim^*_Q y$), then $\iota(x) \preceq \iota(y)$ (resp.\ $\iota(x)\preceq^*\iota(y)$), as we wanted.
    
    The proof that $\iota$ is order-reflecting is similar. Given $x,y\in \dot V(Q)$ and $a$ such that $x,y\in a$ and $a=\mathsf{TC}(a)$, we prove simultaneously by induction that
    \[
        \forall \alpha (\forall x',y'\in a (x',y'\in \dot V_\alpha (Q) \wedge \iota(x') \preceq \iota(y') \rightarrow x'\lesssim_Q y')) 
    \]
    and
    \[
        \forall \alpha (\forall x',y'\in a (x',y'\in \dot V_\alpha (Q) \wedge \iota(x') \preceq^* \iota(y') \rightarrow x'\lesssim^*_Q y')). 
    \]
    We suppose the claim holds for $\beta<\alpha$, and we prove it for $\alpha$.
    Again, we proceed by cases. 
    \begin{itemize}
        \item If $x',y'\in Q$, the claim is obvious.
        \item If $x'\in Q$, $y'\not\in Q$, let $f: \iota(x') \hookrightarrow \iota(y')$ be the (weak) embedding that we are assuming to exist, and let $n\in\omega$ be minimal such that $f(0)\in \sum_{i\leq n} | \iota(\mathsf{enu}(y',h(i)))|$. Let $y''=\mathsf{enu}(y',h(n))$, then $\iota(x')\preceq \iota(y'')$ (equivalently, $\iota(x')\preceq^* \iota(y'')$), and we conclude by inductive assumption.
        \item If $x',y'\not\in Q$, notice that the assumption $\iota(x')\preceq \iota(y')$ (resp.\ $\iota(x')\preceq^* \iota(y')$) entails that for every $x'' \in x'$ there is an $i_{x''}\in\omega$ such that a tail of $\iota(x'')$ (weakly) embeds into $\iota(y',h(i_{x''}))$. Since we know that $\iota(x'')$ is indecomposable, it follows that $\iota(x'')\preceq \iota(y',h(i_{x''}))$ (resp.\ $\iota(x'')\preceq^* \iota(y',h(i_{x''}))$). By inductive assumption, we conclude that $\forall x''\in x'\exists y''\in y'(x''\lesssim^{(*)}_Q y'')$, and so $x'\lesssim^{(*)}_Q y'$.
    \end{itemize}
    For the weak case, we also have to consider the following.
    \begin{itemize}
        \item If $x'\not\in Q$ and $y'\in Q$, the weak embedding $f: \iota(x') \hookrightarrow \iota(y')$ is the constant $0$ map, and such a map also witnesses that $\iota(x'')\preceq^*\iota(y')$ for every $x''\in x'$. Arguing as we did in \Cref{lem:unwinding}, we conclude that $x'\lesssim^*_Q y'$.
    \end{itemize}
    This concludes the proof that $\iota$ is order-reflecting.
\end{proof}

We now define two maps (one for the non-weak and one for the weak case) going in the other direction, i.e.\ mapping transfinite sequences to elements of $\dot V(Q)$.

\begin{Definition}\label{def:eta}
    $(\prsou)$
  Let $\eta$ ($\eta^*$) be a partial class function mapping each $\seq {\omega^{1+\alpha}} Q$ ($\seqalt {\omega^{1+\alpha}} Q$) to $\dot V_\alpha (Q)$ ($\idcsalt \alpha Q$):
  \begin{enumerate}
  \item $\eta^{(*)}(v)=v(|v|-1)$ if the length $|v|$ is a successor ordinal;
  \item $\eta^{(*)}(v)=\min_{\lesssim^{(*)}_Q} \{ \{\eta(v\rst_{[\alpha,\beta)})\mid \alpha<\beta<|v|\} \mid \alpha<|v|\}$, otherwise.
  \end{enumerate}
  Here  $\eta^{(*)}(v)$ will be undefined if either one of values $\eta(v\rst_{[\alpha,\beta)})$ is undefined or the minimum does not exist.
\end{Definition}
Formally in $\prsou$ we define $\eta^{(*)}$ as a total function that simply takes value $\emptyset$ on the arguments where it should not be defined. We compute $\eta(v)$ by computing the $\eta \rst \{v\rst_{[\alpha,\alpha+\beta)}\mid \alpha<|v| \text{ and } \alpha+\beta\le |v|\}$ by recursion on $\beta\le |v|$. 

Intuitively, the maps $\eta$ and $\eta^*$ build an element of $\dot V(Q)$ the correspond to the elements of $Q$ that appear cofinally often in a sequence. This intuition is supported by the following Lemma.

For the sake of readability, we define the following relation: we write $v\sqsubseteq_s u$ if $v$ is a subsegment of $u$, namely if there exist $\alpha,\beta\leq |u|$ such that $=u\rst_{[\alpha,\beta)}$, and $v\sqsubseteq_t u$ if $v$ is a tail of $u$, namely if there exists $\alpha<|u|$ such that $v=u\rst_{[\alpha,|u|)}$. 

\begin{Lemma}\label{lem:useful-tails}
    Let $u$ be a sequence on which $\eta^{(*)}$ is defined, and let $v$ be a tail of $u$. Then $\eta^{(*)}$ is defined on $v$ and moreover $\eta^{(*)}(u) \sim^{(*)}_Q\eta^{(*)}(v)$.
\end{Lemma}
\begin{proof}
    Let us fix $u$, we prove the following by induction:
    \[
        \forall \alpha (\forall u',v'\sqsubseteq_s u (|u'|\leq \alpha \wedge v'\sqsubseteq_t u' \wedge u'\in\dom(\eta) \rightarrow \eta(u')\sim_Q\eta(v'))).
    \]
    Suppose the claim holds for every $\beta<\alpha$, we prove it for $\alpha$. 

    Let $A=\{ \{\eta(u'\rst_{[\gamma,\delta)}) : \gamma< \delta< |u'|\} : \gamma< |u'|\}$ and $B=\{ \{\eta(v'\rst_{[\gamma,\delta)}) : \gamma< \delta< |v'|\} : \gamma< |v'|\}$, so that $\eta(u')=\min_{\lesssim_Q} A$ and $\eta(v')=\min_{\lesssim_Q}B$. Notice that, since $v'\sqsubseteq_t u'$, $B\subseteq A$, so we just have to show that for every element $x\in A$ there is $y\in B$ with $y\lesssim_Q x$. Even more is true: for every $x\in A$, either $x\in B$, or $y\lesssim_Q x$ for every $y\in B$. To see this, fix $y\in B$, suppose that $x\not \in B$, and say $x=\{\eta(u'\rst_{[\gamma,\delta)}) : \gamma< \delta< |u'|\}$ for a certain $\gamma$. Then, for every sufficiently large $\delta$, $u'\rst_{[\gamma,\delta)}$ has a tail $v''$ such that $\eta(v'')\in y$, and every element of $y$ can be obtained this way. We conclude by inductive assumption that $\eta(v'')\sim_Q\eta(u'\rst_{[\gamma,\delta)})$, and so that $y\lesssim_Q x$.

    The case of $\eta^*$ is identical.
\end{proof}
Notice that the Lemma above entails that if $\eta^{(*)}(u)$ is defined, then $\eta^{(*)}$ is also defined on every subsegment of $u$.

Ideally, we would like the maps $\eta$ and $\eta^*$ to enjoy the same strong properties as $\iota$, namely order-preservation and order-reflection. Unfortunately, this is quite clearly not the case with the usual comparison of sequences: consider as urelements the qo $Q=(\{0,1\},=)$. Then, the sequences $01$ and $11$ are both mapped by $\eta$ to the urelement $1$, but $01\not\preceq 11$, showing a failure of order-reflection. Similarly, for the sequences $0$ and $01$ we have $0\preceq 01$, yet $\eta(0)=0\not\lesssim_Q 1 =\eta(01)$, hence order-preservation fails as well. Yet, $\eta$ and $\eta^*$ do have a nice behavior provided we change the notion of comparison of sequences.

\begin{Definition}
($\prsou$)
  For transfinite sequences $u,v$ over a qo $Q$ we say that $u$ \emph{(weakly) cofinally embeds} into $v$ if for every tail $v'$ of $v$ there is a tail $u'$ of $u$ such that $u'\preceq v'$. We denote this relation as $u\cofemb v$ ($u\cofembw v$).
\end{Definition}
It is easy to see that this new relation is reflexive and transitive. 

It is interesting to see what happens if we restrict these new orders to the indecomposable sequences.

\begin{Lemma}\label{lem:cofemb-indec}
($\prsou$)
For every sequences $u$ and $v$, the following hold:
\begin{enumerate}
    \item if $u$ is (weakly) indecomposable and $u\cofemb^{(*)} v$, then $u\preceq^{(*)} v$, and
    \item if $v$ is (weakly) indecomposable and $u\preceq^{(*)} v$, then $u\cofemb^{(*)} v$.
\end{enumerate}

Hence, on $\indseq {\omega^{1+\alpha}} Q$ ($\indseqalt {\omega^{1+\alpha}} Q$) the relations $\preceq$ ($\preceq^*$) and $\cofemb $ ($\cofembw $) coincide.
\end{Lemma}
\begin{proof}
   $u\cofemb^{(*)} v$ implies that a tail $u'$ of $u$ (weakly) embeds into $v$. By (weak) indecomposabilty of $u$, $u\preceq^{(*)} u'$, and so $u\preceq^{(*)} v$. 
    
   By (weak) indecomposability, $v$ (weakly) embeds into all of its tails, and hence so does $u$. This proves that $u\cofemb^{(*)} v$, where  the tail of $u$ we choose for every tail of $v$ is $u$ itself.
\end{proof}

A good way to establish some fundamental properties of $\eta^{(*)}$ is to study its relationship with $\iota$.

\begin{Lemma}\label{lem:iota-eta}
    $(\prsou(\mathsf{enu}))$ 
    For every $x\in \dot V(Q)$ and every sequence $u$ on which $\eta^{(*)}$ is defined, the following hold:
    \begin{enumerate}
        \item\label{eq-i-e-1} $x\lesssim^{(*)}_Q \eta^{(*)}(u) \rightarrow \iota(x) \cofemb^{(*)} u$, and
        \item\label{eq-i-e-2} $\eta^{(*)}(u) \lesssim^{(*)}_Q x \rightarrow u \cofemb^{(*)} \iota(x)$.
    \end{enumerate}
\end{Lemma}
\begin{proof}
    Let us fix $x\in \dot V(Q)$ and a sequence $u$ on which $\eta^{(*)}$ is defined. Let $a$ be a set such that $x\in a$ and $a=\mathsf{TC}(a)$. 
    
    We start proving \Cref{eq-i-e-1}. We prove the following simultaneously by induction.
    \[
        \forall \alpha (\forall u'\sqsubseteq_s u \forall x'\in a(
        |u'|\leq\alpha \wedge x'\lesssim_Q \eta(u') \rightarrow \iota(x') \cofemb u')),  
    \]
    and 
    \[
        \forall \alpha (\forall u'\sqsubseteq_s u \forall x'\in a(
        |u'|\leq\alpha \wedge x'\lesssim^*_Q \eta(u') \rightarrow \iota(x') \cofemb^* u')).
    \]
    Suppose that the claim holds for all $\beta<\alpha$, we prove it for $\alpha$.

    We focus on the non-weak case first, and consider different cases.
    \begin{itemize}
        \item If $x',\eta(u')\in Q$, then $u'$ is a sequence of successor length and $u'(|u'|-1)=p\in Q$ with $q\leq_Q p$, and hence $\iota(x')\cofemb u$ obviously holds.
        \item Suppose $x'\in Q$ and $\eta(u')\not\in Q$. Let $\gamma_0<|u'|$ be such that $\eta(u')=\{\eta(u'\rst_{[\gamma_0,\delta)}):\gamma_0<\delta<|u'|\}$: notice that for every $\gamma>\gamma_0$, 
        \[
            A:=\{\eta(u'\rst_{[\gamma_0,\delta)}):\gamma_0<\delta<|u'|\} \sim_Q
            B:=\{\eta(u'\rst_{[\gamma,\delta)}):\gamma<\delta<|u'|\}.
        \]
        Indeed, by \Cref{lem:useful-tails}, $B\lesssim_Q A$, and the minimality of $A$ forces $A\sim_Q B$.  Hence, for every $\gamma>\gamma_0$ (and hence for every $\gamma<|u'|$), we can find a $\delta>\gamma$ such that $x\lesssim_Q \eta(u'\rst_{[\gamma,\delta)}))$, so by inductive assumption $(x')=\iota(x')\cofemb u'\rst_{[\gamma,\delta)})$, which means $\iota(x')\preceq u'\rst_{[\gamma,\delta)})$. By letting $\gamma$ vary, this yields that $\iota(x')\cofemb u'$.
        \item Suppose then that $x'$ and $\eta(u')$ are both sets. This argument is a slightly more complicated version of the previous one. 
        Let $\gamma_0<|u'|$ be such that $\eta(u')=\{\eta(u'\rst_{[\gamma_0,\delta)}):\gamma_0<\delta<|u'|\}$.
        The assumption means that for every $x''\in x'$ there is $\eta(u'\rst_{[\gamma_0,\delta)})\in \eta(u')$ such that $x''\lesssim_Q \eta(u'\rst_{[\gamma,\delta)})$. Arguing as in the previous case, we can conclude that for every $\gamma<|u'|$ there is $\delta<|u'|$ with $x''\lesssim_Q \eta(u'\rst_{[\gamma,\delta)})$, for which the inductive assumption gives $\iota(x'')\cofemb u'\rst_{[\gamma,\delta)}$ and \Cref{lem:cofemb-indec} gives $\iota(x'')\preceq u'\rst_{[\gamma,\delta)}$. Given any tail $u''$ of $u'$, it is then immediate to build an embedding of $\iota(x')$ into $u''$: since $\iota(x')$ is a concatenation with repetitions of the $iota(x'')$, we can just add the embeddings of $\iota(x'')$ into $u'\rst_{[\gamma,\delta)}$, for some appropriately large $\gamma$.
    \end{itemize}
    The case of $\eta^*$ is almost identical, we only need to consider the additional case
    \begin{itemize}
        \item If $\eta(u')\in Q$, $x\not\in Q$, we again have that $|u'|$ is a successor,and we conclude by inductive assumption that $\iota(x')\preceq^* u'(|u'|-1)$.
    \end{itemize}
    We move to the proof of \Cref{eq-i-e-2}. We prove the following simultaneously by induction.
    \[
        \forall \alpha (\forall u'\sqsubseteq_s u \forall x'\in a(
        |u'|\leq\alpha \wedge \eta(u')\lesssim_Q x' \rightarrow u' \cofemb \iota(x'))),  
    \]
    and 
    \[
        \forall \alpha (\forall u'\sqsubseteq_s u \forall x'\in a(
        |u'|\leq\alpha \wedge \eta(u')\lesssim^*_Q x' \rightarrow u' \cofemb^* \iota(x'))).
    \]
    Suppose that the claim holds for all $\beta<\alpha$, we prove it for $\alpha$.

    We start again with the non-weak case, and consider cases.
    \begin{itemize}
        \item If $\eta(u'),x'\in Q$, then $|u'|$ is a successor and $u'\preceq \iota(x')$ since $u'(|u'|-1)\leq_Q x$. 
        \item If $\eta(u')\in Q$ and $x'\not\in Q$, $|u'|$ is again a successor, $u'(|u'|-1)\lesssim_Q x''$ for some $x''\in x'$ and so $u'\cofemb \iota(x')$ since $\iota(x'')$ appears cofinally often in $\iota(x')$.
        \item The interesting case is again when $x',\eta(u')$ are not in $Q$. As before, let $\gamma_0<|u'|$ be such that $\eta(u')=\{\eta(u'\rst_{[\gamma_0,\delta)}):\gamma_0<\delta<|u'|\}$. Let $u'':=u'\rst_{[\gamma_0,|u'|)}$: we claim that $u''\preceq \iota(x')$. If we do this, by indecomposability of $\iota(x')$, we have that $u'\cofemb \iota(x')$, as we wanted. To show this, we start building the "canonical" embedding $f:u''\to \iota(x')$, namely by recursion we map any $\zeta\in \dom u''$ to the smallest $\xi\in \dom \iota(x')$ such that $\xi> f(\zeta')$ for $\zeta'<\zeta$, and $u''(\zeta)\leq_Q \iota(x')(\xi)$. We want to prove that this procedure covers all of $\dom u''$. Suppose it does not, then there is a minimal $\epsilon< |u''|$ such that $f(\epsilon)$ is undefined, i.e.\ there is no $\xi\in \dom \iota(x')$ such that $\xi> f(\zeta)$ for all $\zeta<\epsilon$ and $u''(\epsilon\leq_Q \iota(x')(\xi)$. We now prove that this is a contradiction.

        To do this, we notice the following fact: reasoning as for \Cref{eq-i-e-1}, we see that for every $\gamma<\delta<|u''|$, there is $x''\in x'$ such that $\eta(u''\rst_{[\gamma,\delta)}))\lesssim_Q x''$, which yields $u''\rst_{[\gamma,\delta)}\cofemb \iota(x'')$ by inductive assumption. We conclude that there is $n\in\omega$ such that $u''\rst_{[0,\epsilon)}\cofemb \iota(\mathsf{enu}(x',h(n)))$, and hence there is $\nu<\epsilon$ such that $u''\rst_{[\nu,\epsilon)}\cofemb \iota(\mathsf{enu}(x',h(n)))$. There is $m\in \omega$ such that $f$ witnesses $u''\rst_{[0,\nu)}\preceq \sum_{i\leq m}\iota(\mathsf{enu}(x',h(i)))$, and by the properties of $h$ there is $\ell>m$ such that $h(n)=h(\ell)$. By trivial composition of embeddings, it follows that $u''\rst_{[0,\epsilon)}\preceq \sum_{i\leq\ell}\iota(\mathsf{enu}(x',h(i)))$. It is now sufficient to notice that $u''\rst_{[\epsilon,\epsilon+1)}\preceq \sum_{i>\ell}\iota(\mathsf{enu}(x',h(i)))$ (again by the same fact we stated at the start of the paragraph) to reach the desired contradiction.
        \end{itemize}
        Again, the case of $\eta^*$ is essentially identical, but we have to consider one more case.
        \begin{itemize}
            \item If $\eta^*(u')\not\in Q$ and $x'\in Q$, let $\gamma_0<|u'|$ be such that $\eta(u')=\{\eta(u'\rst_{[\gamma_0,\delta)}):\gamma_0<\delta<|u'|\}$. Then, for every $\gamma\geq\gamma_0$, $\eta(u'\rst_{[\gamma,\gamma+1)})\lesssim^*_Q x'$, which by inductive assumption entails $u'\rst_{[\gamma,\gamma+1)}\preceq^* \iota(x')$: then, $u\rst_{[\gamma_0,|u'|)}\preceq^* x'$, thus proving that $u'\cofembw x'$ 
        \end{itemize}
    With this, the proof is complete.
\end{proof}

The relationship between $\eta^{(*)}$ and $\iota$ can be summarized as follows.

\begin{Lemma}\label{lem:eta-iota-equiv}
    ($\prsou(\mathsf{enu})$)
    The qo's $(\dom(\eta^{(*)}), \cofemb^{(*)})$ and $(\dot V(Q),\lesssim^{(*)})$ are equivalent as categories, and this is witnessed by the monotone maps $\eta^{(*)}$ and $\iota$.

    Equivalently, the following relations hold for every $x\in \dot V(Q)$ and sequence $u\in \dom(\eta^{(*)})$:
    \[
        \iota (\eta^{(*)} (u)) \cofeq^{(*)} u \quad \text{ and } \quad
        \eta^{(*)}(\iota(x)) \sim^{(*)}_Q x.
    \]
\end{Lemma}
\begin{proof}
    Consider $u\in \dom \eta^{(*)}$. By \Cref{lem:iota-eta}, from $\eta^{(*)} \lesssim^{(*)}_Q \eta^{(*)}(u)$ we obtain both $\iota (\eta^{(*)} (u))\cofemb^{(*)} u$ and $u \cofemb^{(*)} \iota (\eta^{(*)} (u))$, from which $\iota (\eta^{(*)} (u)) \cofeq^{(*)} u$ follows immediately. 

    For the other result: by what we just proved, $\iota(x)\cofeq^{(*)} \iota(x)$ implies 
    $\iota (\eta^{(*)} (\iota(x))) \cofeq^{(*)} \iota(x)$. Since both sides are indecomposable sequences, it follows that $\iota (\eta^{(*)} (\iota(x))) \equiv^{(*)} \iota(x)$ (where $u\equiv^{(*)} v$ iff $u\preceq^{(*)} v$ and $v\preceq^{(*)} u$.
    Since we proved that $\iota$ is order-preserving and order-reflecting (with respect to $\preceq^{(*)}$), it follows that $\eta^{(*)}(\iota(x)) \sim^{(*)}_Q x$.
\end{proof}

We now prove the desired result about the order-theoretic properties of the maps $\eta$ and $\eta^*$.

\begin{Lemma}\label{lem:eta-pres-refl}
    ($\prsou(\mathsf{enu})$)
  For any $\alpha$, the map 
  \[\eta^{(*)} {\upharpoonright}_{\seq {\omega^{1+\alpha}} Q}\colon (\seq {\omega^{1+\alpha}} Q \cap \dom(\eta^{(*)}),\cofemb^{(*)})\to \dot V^{(*)}_\alpha (Q)\] 
  is order-preserving and order-reflecting.
\end{Lemma}
\begin{proof}
    We start showing that the range of the functions is correct, namely that elements of $\dot V_\alpha(Q)$ are mapped to elements of $\seq {\omega^{1+\alpha}} Q$, if the map $\eta$ (or $\eta^*$) is defined. 
    We will proceed by induction on the length of the subsegments of a fixed sequence $u$. 

    So fix a sequence $u$ on which $\eta$ is defined, we prove by induction that
    \[
        \forall \alpha (\forall u'\sqsubseteq_s u (|u'|< \omega^{1+\alpha} \rightarrow \eta^{(*)}(u')\in \dot V_\alpha(Q))). 
    \]
    Suppose we proved the claim for $\beta<\alpha$, we prove it for $\alpha$. If $\alpha$ is $0$ or limit, the claim is immediate. If instead $\alpha=\alpha'+1$, $u'$ has a tail $u''$ such that $|u''|\leq\omega^{1+\alpha'}$, and by \Cref{lem:useful-tails} $\eta^{(*)}(u')=\eta^{(*)}(u'')$. By definition of $\eta^{(*)}$ and inductive assumption, $\eta^{(*)}(u'')$ is a set of elements of $\dot V_{\alpha'}(Q)$, and is thus an element of $\dot V_{\alpha}(Q)$.

    Order-preservation and order-reflection then follow from the previous results: given $u,v\in \seq {\omega^{1+\alpha}} Q \cap \dom (\eta^{(*)})$, we know from above that $u\preceq^{(*)} v$ holds if and only if $\iota(\eta^{(*)}(u)) \preceq^{(*)} \iota(\eta^{(*)}(v))$ does, and the latter holds if and only if $\eta^{(*)}(u) \lesssim^{(*)}_Q \eta^{(*)}(v)$ does, by order-preservation and order-reflection of $\iota$. 
\end{proof}

We can be rather precise about the transfinite sequences in the domain of $\eta$: we describe them in the following result. 

\begin{Lemma}\label{lem:dom-eta-indec}
($\prsou(\mathsf{enu})$)
The domain of $\eta^{(*)}$ is precisely the class of sequences $u$ such that for any $0\le \alpha<\beta\le |u|$, $u\rst_{[\alpha,\beta)}$ splits into sum of finitely many (weak) indecomposable sequences.
\end{Lemma}
\begin{proof}
    We start showing that if $u$ has a (weak) indecomposable tail and $\eta^{(*)}$ is defined on every shorter subsegment of $u$, then $\eta^{(*)}$ is defined on $u$. It is clear that if we do this, then every sequence such that every subsegment can be split into finitely many (weak) indecomposable sequences is in the domain of $\eta^{(*)}$. 

    So suppose that $u$ satisfies these assumptions. Since $\eta^{(*)}$ is defined on every subsegment of $u$, we only have to show that the minimum exists. Let $\gamma_0<|u|$ be such that $u\rst_{[\gamma_0,|u|)}$ is indecomposable, we claim that 
    \[
    \{\eta^{(*)}(u\rst_{[\gamma_0,\delta)}) : \gamma_0<\delta<|u|\} \lesssim^{(*)}_Q \{\eta^{(*)}(u\rst_{[\gamma,\delta)}) : \gamma<\delta<|u|\}
    \]
    for every $\gamma_0<\gamma<|u|$: if we do this, $\{\eta^{(*)}(u\rst_{[\gamma_0,\delta)}) : \gamma_0<\delta<|u|\}$ is the minimum we are looking for. So fix a $\gamma$ as above, by \Cref{lem:eta-pres-refl} it suffices to show that for every $\gamma_0<\delta_0<|u|$ there is $\gamma<\delta<|u|$ such that $u\rst_{[\gamma_0,\delta_0)} \cofemb^{(*)} u\rst_{[\gamma,\delta)}$. By (weak) indecomposability of $u\rst_{[\gamma_0,|u|)}$, there is a (weak) embedding of $u\rst_{[\gamma_0,|u|)}$ into $u\rst_{[\gamma,|u|)}$, say that it is witnessed by the function $f: [\gamma_0,|u|)\to [\gamma,|u|)$ between the two domains. Then $u\rst_{[\gamma_0,\delta_0)} \preceq^{(*)} u\rst_{[\gamma,f(\delta_0))}$, and since we can choose this (weak) embedding to be cofinal, we also have $u\rst_{[\gamma_0,\delta_0)} \cofemb^{(*)} u\rst_{[\gamma,f(\delta_0))}$.

    In the other direction, suppose that $\eta^{(*)}(u)$ is defined: from \Cref{lem:eta-iota-equiv}, it follows that $u$ has a (weakly) indecomposable tail, and so in particular $u$ can be decomposed as sum of finitely many (weakly) indecomposable sequences: since, as we have already observed, $\eta^{(*)}(u)$ defined implies that $\eta^{(*)}$ is defined on every subsegment of $u$ as well, it follows that every subsegment of $u$ can be decomposed as finite sum of (weakly) indecomposable sequences. 
\end{proof}



We can also give conditions on $Q$ that make $\eta^{(*)}$ total.

\begin{Lemma}\label{lem:dom-eta}
($\prsou(\mathsf{enu})$)
If $\dot V_\alpha(Q)$ ($\dot V_\alpha^*(Q)$) is well-founded, then $\eta^{(*)}$ is total on $\seq {\omega^{1+\alpha}} Q$ ($\seq {\omega^{1+\alpha}} Q$).
\end{Lemma}
\begin{proof}
     Let us fix $u\in \seq {\omega^{1+\alpha}} Q$, we prove that if $u'\in \dom(\eta^{(*)})$ for every $u'\sqsubset_s u$ of smaller length, then $u\in \dom (\eta^{(*)})$, if we do this the Lemma clearly follows. Since we know that $\eta^{(*)}(u\rst_{[\gamma,\delta)})$ exists for all $0\leq \gamma< \delta <|u|$, we only have to show that the minimum exists. We have seen in \Cref{lem:eta-pres-refl} that the range of $\eta^{(*)}\rst_{\seq {\omega^{1+\alpha}} Q}$ is $\dot V_\alpha(Q)$, so $\{ \{\eta(v\rst_{[\gamma,\delta)})\mid \gamma<\delta<|v|\} \mid \gamma<|v|\}$ is a descending sequence in $\dot V_\alpha(Q)$ ($\dot V_\alpha^*(Q)$): since the latter is well-founded by assumption, this sequence has a minimum.
\end{proof}

We put together all the results we have proved so far in the following Theorems.

\begin{Theorem}\label{theo:big-theo-1}
($\prsou(\mathsf{enu})$)
  For every ordinal $\alpha$, the following are equivalent:
  \begin{enumerate}
  \item\label{bigtheo-1} $Q$ is $\alpha$-wqo;
  \item\label{bigtheo-2} $\dot V_\alpha (Q)$ is well-founded;
  \item\label{bigtheo-3}  $\dot V_\beta (Q)$ are wqo for every $\beta<\alpha$ (in the case of $\alpha>0$);
  \item\label{bigtheo-5}  $\idcsalt \beta Q$ are wqo for every $\beta<\alpha$ (in the case of $\alpha>0$);
  \item\label{bigtheo-6} $\indseq {\omega^{1+\alpha}} Q$ is well-founded;
  \item\label{bigtheo-7} $\indseq {\omega^{1+\beta}} Q$ are wqo for $\beta<\alpha$ (in the case of $\alpha>0$);
  \item\label{bigtheo-8} $\seq {\omega^{1+\alpha}} Q$ is well-founded;
  \item\label{bigtheo-9} $\seq {\omega^{1+\beta}} Q$ are wqo for $\beta<\alpha$ (in the case of $\alpha>0$);
  \item\label{bigtheo-11} $\indseqalt {\omega^{1+\beta}} Q$ are wqo for all $\beta<\alpha$ (in the case of $\alpha>0$);
  \item\label{bigtheo-13} $\seqalt {\omega^{1+\beta}} Q$ are wqo for $\beta<\alpha$ (in the case of $\alpha>0$);
  \end{enumerate}
\end{Theorem}

\begin{proof}
The proof proceeds as follows: 
first, we prove that, for every ordinal $\alpha$, the equivalences \ref{bigtheo-1} $\Leftrightarrow$ \ref{bigtheo-3} $\Leftrightarrow$ \ref{bigtheo-5} $\Leftrightarrow$ \ref{bigtheo-7} $\Leftrightarrow$ \ref{bigtheo-9} $\Leftrightarrow$ \ref{bigtheo-11} $\Leftrightarrow$\ref{bigtheo-13} and \ref{bigtheo-2} $\Leftrightarrow$ \ref{bigtheo-6} $\Leftrightarrow$ \ref{bigtheo-8} hold, as does the implication $\ref{bigtheo-2} \Rightarrow \ref{bigtheo-3}$. Then we show that $\ref{bigtheo-3} \Rightarrow \ref{bigtheo-2}$ by distinguishing the cases that $\alpha$ is a successor or a limit ordinal (the equivalences $\ref{bigtheo-1} \Leftrightarrow \ref{bigtheo-2} \Leftrightarrow \ref{bigtheo-6} \leftrightarrow \ref{bigtheo-8}$ in the case of $\alpha=0$ are obvious, and the others do not make sense).


We have already observed in \Cref{subsec:wind-unwind} that $Q$ is $(\alpha+1)$-wqo if and only if $\dot V_\alpha(Q)$ is wqo if and only if $\idcsalt \alpha Q$ is wqo, which gives $\ref{bigtheo-1} \Leftrightarrow \ref{bigtheo-3} \Leftrightarrow \ref{bigtheo-5}$, by unwinding bad sequences and winding bad arrays. $\ref{bigtheo-3} \Rightarrow \ref{bigtheo-9}$ follows from the fact that, by \Cref{lem:dom-eta}, $\eta$ is defined on $\seq {\omega^{1+\beta}} Q$, and it is an order-reflecting map into $\dot V_\beta (Q)$ for every $\beta<\alpha$. The implication \ref{bigtheo-5} $\Rightarrow$ \ref{bigtheo-13} is proved in the same way. The implications \ref{bigtheo-9} $\Rightarrow$ \ref{bigtheo-7} and \ref{bigtheo-13} $\Rightarrow$ \ref{bigtheo-11} are obvious. Finally, to see that \ref{bigtheo-7} $\Rightarrow$ \ref{bigtheo-3}, it is enough to notice that, for every $\beta<\alpha$, $\iota$ is an order-reflecting map from $\dot V_\beta(Q)$ to $\indseq {\omega^{1+\beta}} Q$, so $\indseq {\omega^{1+\beta}} Q$ wqo implies $\dot V_\beta(Q)$ wqo. \ref{bigtheo-11} $\Rightarrow$ \ref{bigtheo-5} is proved in the same way.

\ref{bigtheo-2} $\Rightarrow$ \ref{bigtheo-8} is again proved by the properties of the function $\eta$: it is total on $\seq {\omega^{1+\alpha}} Q$ by \Cref{lem:dom-eta}, and order-reflecting. \ref{bigtheo-8} $\Rightarrow$ \ref{bigtheo-6} is obvious, and \ref{bigtheo-6} $\Rightarrow$ \ref{bigtheo-2} follows from the fact that $\iota$ is order-preserving and order-reflecting, with $\indseq {\omega^{1+\alpha}} Q$ as range. 

The argument to show that \ref{bigtheo-2} $\Rightarrow$ \ref{bigtheo-3} is analogous to what we have already done in \Cref{lem:vdc-wf-skel}: suppose $\dot V_\beta (Q)$ is not a wqo, for some $\beta<\alpha$, then there is a bad sequence $f:\omega\to \dot V_\beta(Q)$. It follows that $g:\omega\to \dot V_{\beta+1} (Q)$ defined as $g(i)= \{f(j): j\geq i\}$ is a descending sequence in $\dot V_{\beta+1}(Q)$, and so in $\dot V_\alpha(Q)$. 

Now suppose that $\alpha=\beta+1$, and suppose for a contradiction that $\dot V_\alpha(Q)$ is not well-founded, as witnessed by the descending sequence $f:\omega\to \dot V_\alpha(Q)$. Let us define $g:\omega\to \dot V_\beta(Q)$ by putting as $g(i)$ any $x\in \dot V_\beta(Q)$ such that $\forall y\in f(i+1) (x\not\lesssim_Q y)$. It is clear that $g$ is a bad sequences, contradicting that $\dot V_\beta(Q)$ is wqo, and thus proving \ref{bigtheo-3} $\Rightarrow$ \ref{bigtheo-2} in the successor case. 

For $\alpha$ limit, we notice that $\dot V_\alpha(Q)$ is well-founded if all $\dot V_\beta(Q)$ for $\beta<\alpha$ are well-founded as well, since any descending sequence in $\dot V_\alpha(Q)$ is actually a descending sequence in some $\dot V_\beta(Q)$ for $\beta<\alpha$, since $x\lesssim_Q y$ and $y\in \dot V_\beta(Q)$ implies $x\in \dot V_\beta(Q)$.
\end{proof}

In view of Lemmas \ref{lem:cofemb-indec}, \ref{lem:eta-iota-equiv} and \ref{lem:dom-eta} if any of conditions from Theorem \ref{theo:big-theo-1} holds, then 
\begin{enumerate}
    \item the quasi-orders $\dot V_\alpha (Q)$ and $\indseq {\omega^{1+\alpha}} Q$ are equivalent as categories with the equivalence given by the maps $\eta$ and $\iota$;
    \item the quasi-orders $\idcsalt \alpha Q$ and $\indseqalt {\omega^{1+\alpha}} Q$ are equivalent as categories with the equivalence given by the maps $\eta^*$ and $\iota$.
\end{enumerate}

The case of weak embeddability is more subtle.

\begin{Theorem}\label{theo:big-theo-2}
    ($\prsou(\mathsf{enu})$)
    For every ordinal $\alpha$, the following are equivalent:
    \begin{enumerate}
        \item\label{bigtheo-4} $\idcsalt \alpha Q$ is well-founded;
        \item\label{bigtheo-10} $\indseqalt {\omega^{1+\alpha}} Q$ is well-founded;
        \item\label{bigtheo-12} $\seqalt {\omega^{1+\alpha}} Q$ is well-founded;
    \end{enumerate}
    and they all imply that $Q$ is $\alpha$-wqo.
    Moreover, if $\alpha$ is not a limit, $Q$ being $\alpha$-wqo implies any of the above.
\end{Theorem}
\begin{proof}
    The proof proceeds exactly as in the previous Theorem: items \ref{bigtheo-4}, \ref{bigtheo-10} and \ref{bigtheo-12} of the statement of this Theorem behave much like \ref{bigtheo-2}, \ref{bigtheo-6} and \ref{bigtheo-8} in the statement of the previous one. The only argument that does not goes through is the last paragraph of the previous proof: if $\alpha$ is limit and $f$ is an infinite descending sequence in $\idcsalt \alpha Q$, we cannot conclude that $f$ is a descending sequence in $\idcsalt \beta Q$ for some $\beta<\alpha$, since it can now happen that $x\lesssim_Q^* y$ even if $x\in \idcsalt \beta Q$, $y\in \idcsalt \gamma Q$ and $\beta>\gamma$.
\end{proof}

We can extract some additional results if we assume that $Q$ is a set.

\begin{Corollary}
    ($\prsou(\mathsf{enu})$)
    For every $\alpha$ the following are equivalent:
    \begin{enumerate}
        \item\label{bigtheo-1set} $Q$ is $\alpha$-wqo and a set;
        \item\label{bigtheo-14} $\dot V^{dc,f}_\alpha(Q)$ is well-founded and forms a set;
        \item\label{bigtheo-15} $\dot V^{dc,f}_\beta(Q)$ are wqo's and form sets, for $\beta<\alpha$ (in the case of $\alpha>0$).
    \end{enumerate}
\end{Corollary}
\begin{proof}
    \ref{bigtheo-1set} $\Rightarrow$ \ref{bigtheo-14} since we have already observed that $Q$ set implies $\dot V^{dc,f}_\alpha(Q)$ set, and well-foundedness follows from the Theorem above, since $\dot V^{dc,f}_\alpha(Q)$ is a suborder of $\dot V_\alpha(Q)$. \ref{bigtheo-14} $\Rightarrow$ \ref{bigtheo-1set} follows from the obvious fact that $\dot V^{dc,f}_\alpha(Q)$ set implies $Q$ set, \Cref{lem:vdc-wf-skel}, and the Theorem above.

    The equivalence \ref{bigtheo-1set} $\Leftrightarrow$ \ref{bigtheo-15} is analogous.
\end{proof}

\section{Sequences with finite range}\label{sec:tr-seq-fin-ran}

In this Section, we show that transfinite sequences with range a finite subset of $Q$ form a wqo, provided $Q$ is a wqo.

Broadly speaking, we will proceed as follows: first, we notice a general fact, namely that every transfinite sequence over $Q$ with finite range is in the domain of $\eta$: thanks to Higman's Lemma, we can thus restrict to consider \emph{indecomposable} sequences with finite range. Given a bad sequence $(u_i)_{i\in\omega}$ of indecomposable sequences of finite range, $(\eta(u_i))_{i\in\omega}$ is then a bad sequence in $\dot V_\alpha(Q)$, for some $\alpha$. Using our careful unwinding procedure from \Cref{lem:unwinding}, we will extract from this a bad sequence in $Q$, thus getting a contradiction.

\begin{Lemma}\label{lem:eta-fin-def}
    ($\prsou(\mathsf{enu})$)
    If $u$ is a transfinite sequence over $Q$ with finite range, then $\eta(u)$ is defined. 
\end{Lemma}
\begin{proof}
    Let us consider $Q'\subseteq Q$: it is immediate to see that we can build $\dot V(Q')$ as we did with $Q$ (the fact that we are dealing with a subclass of the class of urelements is immaterial). Equally, it is easy to see that if we define the order $\lesssim_{Q'}$ on $\dot V(Q')$ as in \Cref{def:lesssim}, we get the restriction of $\lesssim_Q$ to $\dot V(Q')$. 

    Now, the definition of $\eta$ "relativized to $Q'$" is 
      \begin{enumerate}
  \item $\eta_{Q'}(u)=u(|u|-1)$ if the length $|u|$ is a successor ordinal;
  \item $\eta_{Q'}(u)=\min_{\lesssim_{Q'}} \{ \{\eta_{Q'}(u\rst_{[\alpha,\beta)})\mid \alpha<\beta<|u|\} \mid \alpha<|u|\}$, otherwise,
  \end{enumerate}
  for every $u$ with range in $Q'$, so by what we said in the paragraph above, an easy induction shows that for every transfinite sequence $u$ with range in $Q'$,
  \[
    \eta_{Q'}(u)= \eta(u).
  \]
  In our case, $u$ has finite range, so $Q'$ is finite: as we have observed, $\prsou$ proves that $Q'$ is bqo. Then, by \Cref{lem:dom-eta}, $\eta_{Q'}(u)$ is defined, and so is $\eta(u)$.
\end{proof}

\begin{Lemma}\label{lem:eta-fin-comb}
    ($\prsou(\mathsf{enu})$)
    If $Q$ is wqo, then $(\{ x\in \dot V(Q): \supp(x) \text{ is finite}\}, \lesssim_Q)$ is wqo.
\end{Lemma}
\begin{proof}
    Suppose for a contradiction that the claim is false, and let $(x_i)_{i\in\omega}$ be a bad sequence in $\dot V(Q)$ such that all the $x_i$ have finite support. Let $\alpha$ be such that $x_i\in \dot V_\alpha(Q)$ for every $i\in \omega$. \Cref{lem:unwinding} then gives us a front $F$ and a bad array $g: F\to Q\times\omega$ obtained from the bad sequence $f: i\mapsto x_i$.

      Let $F'$ be the set of finite strings defined as 
  \[
    F':= \{ \{n\} \cup \sigma_0\cup \sigma_1\in [\omega]^{<\omega}:
    n<\sigma_0<\sigma_1 \text{ and } \{n\} \cup \sigma_0,\{n \} \cup \sigma_1\in F\}.
  \]
  We verify that $F'$ is a front:
  \begin{itemize}
  \item Clearly, $\bigcup F'=\bigcup F$.
  \item Let $X\in [\bigcup F]^\omega$, and let $\sigma\in F$ be such that
    $\sigma \sqsubseteq X$. If $\sigma=\{ n\}$ for some $n\in\omega$, then
    $\{n\}\in F'$ as well, which proves that $X$ has an initial segment in $F'$.
    Otherwise, $\sigma=\{n \} \cup \sigma_0$, for some non-empty $\sigma_0$.
    Then, it is enough to consider the $\tau\in F$ such that
    $\tau \sqsubseteq X\setminus \sigma_0$ in order to find a string $\sigma_1$
    (namely, $\tau\setminus \{n\}$) such that $\{n\} \cup \sigma_0\cup \sigma_1\in F'$
    and  $\{n\} \cup \sigma_0\cup \sigma_1\sqsubseteq X$.
  \item Finally, suppose for a contradiction that there are $\tau_0,\tau_1\in F'$
    such that $\tau_0\sqsubseteq \tau_1$. Say that
    $\tau_i=\{n^i\}\cup \{\sigma^i_0\} \cup \{ \sigma^i_1\}$ for $i<2$, then by the
    assumptions we have $n^0=n^1$, and so either
    $\sigma^0_0\sqsubseteq \sigma^1_0$ or $\sigma^1_0\sqsubseteq \sigma^0_0$,
    and both case contradict the assumption that $F$ is a front. 
  \end{itemize}

 We let $g': F\to Q$ be the projection of $g$ on its first component. 
  We then define a coloring $c:F'\to 2$ as follows:
  \[
    c(\{n\} \cup \sigma_0 \cup \sigma_1)=
    \begin{cases}
      0 & \text{ if } g'(\{n\} \cup \sigma_0)= g'(\{n\} \cup \sigma_1) \\
      1 & \text{ otherwise.}
    \end{cases}
  \]
  We can then apply Clopen Ramsey Theorem to the coloring $c$ above, obtaining
  an infinite set $H\subseteq \omega$ such that the
  subfront $F'\rst_H$ of $F'$ is $c$-homogeneous. 

  We claim that $F'\rst_H$ is necessarily $c$-homogeneous for color $0$. We prove the claim by
  contradiction: suppose that there is a subfront $G$ of $F'$ that is $c$-homogeneous
  for $1$. We notice that, for every $n\in \bigcup G$, $\{n\}\not\in G$, since
  $c(\{n\})=0$. In particular, let $m$ be the minimal element of $\bigcup G$: then, we
  can find an infinite sequence $\sigma_0<\sigma_1<\dots$ of finite strings such that
  $\{m\}\cup \sigma_i\in G$.
  By definition of $G$, for every $i\neq j\in \omega$, we have $g'(\{m\}\cup\sigma_i)\neq g'(\{m\}\cup \sigma_j)$. But, as noticed in \Cref{remark:supp}, $g(\{m\}\cup \sigma_i)\in \supp(x_m)$ for every $i\in \omega$, and we are assuming that $\supp(x_m)$ is finite. This gives a contradiction and proves that $F'\rst_{H}$ is homogeneous for $0$.
  
  Let $f':H\to Q$ be defined as $f(n)=g'(\{n\}\cup \sigma)$, for any $\sigma$ such that $\{n\}\cup \sigma \in \{ \tau \in F: \tau\subseteq H\}$: it follows from the definition of $H$
  that $f'$ is well-defined. By the assumption that $Q$ is wqo, we can find an infinite set $A\subseteq H$ such that $\ran f'\rst_A$ is well-ordered.

  Finally, consider the restriction $g\rst_A: F\rst_A \to Q\times \omega$. Since $g$ was assumed to be a bad array, then so is every restriction of it. But, by the
  definition of $A$, the range of $g\rst_A$ projected on the first component is well-ordered by $\leq_Q$. Moreover, we have already observed that the product of well-orders is a bqo, which means that $g\rst_A$ is good. This final contradiction proves the claim.
\end{proof}

As outlined above, we just have to put the previous results together.

\begin{Corollary}
    ($\prsou(\mathsf{enu})$)
    Let $(u_i)_{i\in\omega}$ be a sequence of transfinite sequences over $Q$ with finite range. Then, $(u_i)_{i\in\omega}$ is good.
\end{Corollary}
\begin{proof}
    As proved in \Cref{lem:eta-fin-def}, $\eta$ is defined on every $u_i$. Hence, by \Cref{lem:dom-eta-indec}, every $u_i$ is the sum of finitely many indecomposable sequences, say $u_i=\sum_{j\leq n_i} v^i_j$, for some integer $n_i$. By Higman's Lemma, it is sufficient to prove that the set $\{ v^i_j\}$ ordered by embeddability is wqo.
    
    It is easily seen that $\supp(\eta(v^i_j))\subseteq \ran u_i$. Hence, by \Cref{lem:eta-fin-comb}, $(\{\eta(v^i_j))\},\lesssim_Q)$ is wqo, and hence $(\{ v^i_j\},\cofemb)$ is wqo. Since every $v^i_j$ is indecomposable, this is equivalent to $(\{ v^i_j\},\preceq)$ being wqo, as we wanted. 
\end{proof}

\appendix

\section{Interpreting $\prsou$ in $\mathsf{ATR}_0$}\label{app:prosu-in-atr}

In this section we sketch how $\prsou$ is interpreted in $\mathsf{ATR}_0$. Unfortunately, although interpretations of $\mathsf{PRS}\omega$ in $\mathsf{ATR}_0$ have been used before, we haven't found a presentation of how precisely the interpretations of primitive-recursive set functions are defined within $\mathsf{ATR}_0$. The goal of this section is to sketch this construction.

As we mentioned before our default embedding of the language of $\prsou$ into the language of $\mathsf{ATR}_0$ is given by the interpretation where sets are interpreted as points within well-founded extensional binary relations with $\mathbb{N}$ as the urelements. More formally, sets or urelements are interpreted as pairs $(A,a)$, where $A$ is a structure with a well-founded binary relation $\in^A$ and a family of pairwise distinct constants $(c_i^A)_{i<\omega}$ such that $A$ is extensional relative to the interpretation of urelements as the constants $c_i$, i.e.
$$\{z\in A\mid z\in^A x\}=\{z\in A\mid z\in^A y\}\;\Rightarrow\;x=y\text{, for any $x,y\in A\setminus \{c_i^A\mid i\in\mathbb{N}\}$ and}$$
\[x\not\in^A c_i^A\text{, for any $i\in\mathbb{N}$}.\]
An amalgamation of two structures $A,B$ as above is a triple $(C,f_A,f_B)$, where $C$ is a structure with the same properties and the functions $f_A\colon A\to C$ and $f_B\colon B\to C$ are initial embeddings (injective maps $f$ preserving all the constants, preserving and reflecting the binary relation and having the property that for any $x\in\mathsf{ran}(f)$ its whole $\in$-downwards cone is also in $\mathsf{ran}(f)$). The equality $(A,a)=^I(B,b)$ is interpreted as existence of an amalgamation $(C,f_A,f_B)$ such that $f_A(a)=f_B(b)$. The membership $(A,a)\in^I(B,b)$ is interpreted as existence of an amalgamation $(C,f_A,f_B)$ such that $f_A(a)\in^C f_B(b)$. A pair $(A,a)$ represents a urelement if $a$ is one of $c_i^A$.

Formally, to complete the definition of the interpretation we also need to define the interpretations of all the primitive-recursive set functions. Under the name of $\Sigma$-recursion Barwise basically shows clousure of $\Sigma$-functions under primitive recursion within $\mathsf{KPU}$. Unfortunately for us, the proof of this general closure property relies on the principle of $\Sigma$-collection which is not provable for our interpretation of set theory in $\mathsf{ATR}_0$.  However, we will get around this by using some properties specific to the primitive recursive set functions. For any primitive recursive set function $f(\vec{x})$, as we will show below, there is always some number number $k$, which we call the \emph{rank} of the function, such that the natural $\Sigma$-definition of $f(\vec{x})=y$ (the definition could be extracted from Barwise's book \cite[Sections~I.V~and~I.VI]{book-barwise}) gives a total function when relativized to any $L_{\varphi_k(\alpha)}(a)$ for any set $a$ whose rank $\mathsf{rk}(a)<\varphi_k(\alpha)$; here $\varphi_k$ is the $k$-th Veblen function. 

Let us prove by induction on construction of the primitive recursive set function $f(\vec{x})$ that all of them have finite rank. The only non-trivial case are the constructions of the functions by primitive recursion. For a function $f(y,\vec{z})=R[t(x,y,\vec{z})](y,\vec{z})$, where all functions in $t$ have rank $\le k$ we are claiming that $f$ is of the rank $\le k+1$. Indeed, let us track the computation of $f(y,\vec{p})$ within some $L_{\varphi_{k+1}(\alpha)}(a)$, where $\mathsf{rk}(a)<\varphi_{k+1}(\alpha)$ and $\vec{p}\in L_{\varphi_{k+1}(\alpha)}(a)$. For this we fix a supplementary ordinal $\beta<\varphi_{k+1}(\alpha)$ that is bigger than the rank of $a$ and all components of $\vec{p}$. For this $\beta$ we, of course, have the property that for any $\gamma<\varphi_{k+1}(\alpha)$ the computations of $t(x,y,\vec{p})$ could be performed within $L_{\varphi_k(\beta+\gamma)}(a)$ (and since $\varphi_{k+1}$ enumerates fixed points of $\varphi_k$, we have $\varphi_k(\beta+\gamma)<\varphi_{k+1}(\alpha)$). Notice that now by $\in$-induction on sets $b\in L_{\varphi_{k+1}(\alpha)}(a)$ we easily show that the computation of $f(b,\vec{p})$ could be performed within $L_{\varphi_k(\beta+\mathsf{rk}(b)+1)}(a)$. 

This facts allows us to define interpretation of primitive-recursive set functions within $\mathsf{ATR}_0$ and furthermore verify that our interpretation indeed satisfies the axiom scheme $\Delta_0$-Foundation. First note that $\mathsf{ATR}_0$ proves that for any well-order $\alpha$, the naturally defined linear orders $\varphi_k(\alpha)$ are also well-ordered (a proof is given in \cite{marcone2011veblen}, though probably the fact was known since the introduction of $\mathsf{ATR}_0$). Furthermore, given any code of a set $(A,a)$ clearly we could construct in $\mathsf{ATR}_0$ a well-order $\alpha$ equal to the rank of $a$ in $A$. And next using arithmetical transfinite recursion along $\varphi_k(\alpha+1)$ we could construct $(B,b)$ representing $L_{\varphi_k(\alpha)}(a)$. We define the interpretation of $f$ of the rank $k$ on $\in^A$-elements of $a$ according to the relativization of the natural $\Sigma$ definition of $f$ to $b$ within the structure $B$. If some $\Delta_0$-formula involves only functions of ranks $\le k$, then its truth value on $\in^A$-elements of $a$ will be arithmetical relative to $(B,b)$, which will allow us to verify $\Delta_0$-Foundation.

\section{Forcing Global Enumerating Function}\label{app:forcing-enu}

Our goal in the section will be to prove the following conservation result
\begin{Theorem}\label{cons_enu}
    The set of consequences of the theory $\prsou(\mathsf{enu})$ in the language of $\prsou$ is precisely $\prsou+\text{``every set is countable''}$.
\end{Theorem}

The idea of the proof is to define forcing within $\prsou+\text{``every set is countable''}$ that would force the function $\mathsf{enu}$.

In $\prsou$ we define a \emph{partial enumerating function} $e$ to be a (set) function whose domain consists of non-empty sets and that maps each $x$ in $\mathsf{dom}(e)$ to a surjective function $e(x)\colon \omega \to x$. We order partial enumerating functions by the expansion order $\supseteq$.

For every $\prsou(\mathsf{enu})$-term $t(\vec{x})$ we define three $\prsou$-terms:
\begin{enumerate}
    \item $t^d(e,\vec{x})$ that returns $1$ if the partial enumerating function $e$ has enough information to compute $t(\vec{x})$;
    \item $t^p(e,\vec{x})$ that computes $t(\vec{x})$ using $e$.
\end{enumerate}
Formally, we make the definition by recursion on the construction of  $\prsou(\mathsf{enu})$-terms:
\begin{enumerate}
    \item for basic functions $f(\vec{x})$ other than $\mathsf{enu}$ we put \[f^d(e,\vec{x})=1\;\;\text{and}\;\;f^p(e,\vec{x})=f(\vec{x}).\]
    \item for $\mathsf{enu}(x,y)$:
    \[\mathsf{enu}^d(e,x,y)=\begin{cases}1 &\text{if either $x\in\mathsf{dom}(e)$, or $x\in U$, or $x=\emptyset$}\\0 &\text{otherwise}\end{cases}\]
    \[\mathsf{enu}^p(e,x,y)=\begin{cases}e(x)(y) &\text{if $x\in\mathsf{dom}(e)$ and $y\in \omega$}\\0 &\text{otherwise}\end{cases}\]
    \item for terms $t(\vec{x})$ of the form $f(\vec{u}(\vec{x}))$ we put 
    \[t^d(e,\vec{x})=f^d(e,\vec{u}^p(e,\vec{x}))\prod_{u\text{ in }\vec{u}}  u^d(e,\vec{x})\;\;\text{and}\;\;t^p(e,\vec{x})=f^p(e,\vec{u}^p(\vec{x})).\]
    \item for a term $t(y,\vec{z})$ of the form $R[u(x,y,\vec{z})](y,\vec{z})$ we put 
    \[t^d(e,y,\vec{z})=u^d(e,\bigcup\{t^p(e,y',\vec{z})\mid y'\in y\},y,\vec{z})\prod_{y'\in y}t^d(e,y',\vec{z}),\]
    \[t^p(e,y,\vec{z})=u^p(e,\bigcup\{t^p(e,y',\vec{z})\mid y'\in y\},y,\vec{z}).\]
\end{enumerate} 


By recursion on construction of formulas $\varphi(\vec{x})$ of the language of $\prsou(\mathsf{enu})$ we define formulas $e\Vdash \varphi(\vec{x})$ of the language of $\prsou$, where $e$ is an additional free variable, whose intended range are partial enumerating functions:
\begin{enumerate}
    \item $e\Vdash t(\vec{x})=u(\vec{x})$ is $\forall e'\supseteq e\exists e''\supseteq e'(t^d(\vec{x})=u^d(\vec{x})=1\land t^p(\vec{x})=u^p(\vec{x}))$;
    \item $e\Vdash t(\vec{x})\in u(\vec{x})$ is $\forall e'\supseteq e\exists e''\supseteq e'(t^d(\vec{x})=u^d(\vec{x})=1\land t^p(\vec{x})\in u^p(\vec{x}))$;
    \item $e\Vdash \varphi\land \psi$ is $(e\Vdash \varphi)\land (e\Vdash \psi)$;
    \item $e\Vdash \forall x\;\varphi(x)$ is $\forall x(e\Vdash \varphi(x))$;
    \item $e\Vdash \lnot\varphi$ is $\forall e'\supseteq e\;\lnot (e\Vdash \varphi)$.
\end{enumerate}
This definition is simply an instance of forcing interpretation \cite{pakhomov2022finitely,avigad2004forcing} and thus forceability is monotone, i.e.\
\[\prsou\vdash e'\supseteq e \land (e\Vdash \varphi)\to (e'\Vdash \varphi),\]
forceability of each formula by $e$ is equivalent to the density of $e'\supseteq e$ that force this formula, i.e.\
\[\prsou\vdash (e\Vdash \varphi)\mathrel{\leftrightarrow} \forall e'\supseteq e\exists e''\supseteq e'(e''\Vdash \varphi),\]
and forceability is closed under classical first-order deduction, i.e. if $\varphi\to \psi$ is provable in classical first-order logic, then
\[\prsou\vdash (e\Vdash \varphi)\to (e'\Vdash \psi).\]

Too prove Theorem \ref{cons_enu} we will need the following two lemmas:

 \begin{Lemma}\label{forcing_conservativity}
    For any $\prsou$-formula $\varphi(\vec{x})$, the theory $\prsou$ proves 
    $$\varphi(\vec{x})\mathrel{\leftrightarrow} (\emptyset\Vdash \varphi(\vec{x})).$$
\end{Lemma}

\begin{Lemma}\label{forcing_axioms}
    For any $\prsou(\mathsf{enu})$-axiom $\varphi$ that is not $\prsou$-axiom, the theory $\prsou+\text{``every set is countable''}$ proves 
    $$\emptyset\Vdash \varphi.$$
\end{Lemma}

Having these two lemmas we reason as follows to prove Theorem \ref{cons_enu}. We consider a theorem $\varphi$ of $\prsou(\mathsf{enu})$ in the language of $\prsou$ and claim that it is already a theorem of $\prsou +\text{``every set is countable''}$. By compactness theorem there are finitely many axioms $\psi_1,\ldots,\psi_n$ such that  $\psi_1\land \ldots \land \psi_n\to \varphi$ is a theorem of first-order logic. Then $\prsou +\text{``every set is countable''}$ proves $\emptyset\Vdash \psi_1\land \ldots \land \psi_n$, where we obtain provability of individual $\emptyset\Vdash \psi_i$ either by Lemma \ref{forcing_conservativity} for axioms of $\prsou$ or  by Lemma \ref{forcing_axioms} for other axioms. By closure of forcing under first-order deduction we get that $\prsou +\text{``every set is countable''}$ proves  $\emptyset\Vdash \varphi$. Thus by Lemma \ref{forcing_conservativity} the theory $\prsou +\text{``every set is countable''}$ proves $\varphi$.

We prove Lemma \ref{forcing_conservativity} by a straightforward induction on construction of formulas. For the rest of the section our goal is to prove Lemma \ref{forcing_axioms}. For this end we will need several supplementary lemmas.

\begin{Lemma} For any $\prsou(\mathsf{enu})$-term $t(\vec{x})$ there is a $\prsou$-term $t^c(e,s)$ such that $\prsou$ proves that if $e$ is a partial enumerating function, $s\in\mathsf{dom}(e)$, $s$ is a transitive set, and $\vec{x}\in s$, then $t^d(t^c(e,s),\vec{x})=1$ and $t^c(e,s)$ is a partial enumerating function extending $e$.\end{Lemma}
\begin{proof}
We define terms $t^c(e,s)$ by induction on construction of $t$. 

For all the basic functions $f(\vec{x})$ we put $f^c(e,s)=e$. 

If $t(\vec{x})$ is $f(u_1(\vec{x}),\ldots,u_n(\vec{x}))$, then we put 
$$t^c(\vec{x})=f^c(e',\mathsf{TC}(\{u_i^p(e',\vec{x})\mid 1\le i\le n, \vec{x}\in s\})),$$
where $e'=u_n^c(\ldots u_2(u_1^c(e,s),s)\ldots,s)$.

If $t(y,\vec{z})$ is $R[u(x,y,\vec{z})](y,\vec{z})$, then we define the function $t^c(e,s)$ as follows. 

First we define function $\mathsf{oemb}(e,s,q)$ such that if $s$ is a transitive set, $e$ is a partial enumerating function with $s\in\mathsf{dom}(e)$, and $q\in s\cup \{s\}$, then $\mathsf{oemb}(e,s,q)$ is a bijective $\in$-preserving map from $\mathsf{TC}(q)$ to an ordinal. We compute $\mathsf{oemb}$ by $\in$-recursion on $q$: If $q$ is $\emptyset$ or a urelement, then $\mathsf{oemb}(e,s,q)=\emptyset$. Otherwise, on the recursion step in order to compute $\mathsf{oemb}(e,s,q)$ we compute sequences of expanding functions $(f_i)_{i<\omega}$ and increasing ordinals $(\alpha_i)_{i<\omega}$:
\begin{itemize}
    \item $f_0=\emptyset$;
    \item $\alpha_i=\mathsf{ran}(f_i)$;
    \item if $e(s)(i)\not\in q\setminus \mathsf{dom}(f_i)$, then $f_{i+1}=f_i$;
    \item if $e(s)(i)=q'\in q\setminus \mathsf{dom}(f_i)$, then we define $f_{i+1}$ as follows:
    \begin{itemize}
    \item we consider a bijective  $\in$-preserving map $g=\mathsf{oemb}(e,s,q')$ with $\mathsf{dom}(g)=\mathsf{TC}(q')$,
    \item we consider $g'=g{\upharpoonright} (\mathsf{TC}(q')\setminus\mathsf{dom}(f_i))$,
    \item we consider $g''=\pi\circ g'$, where $\pi$ is the unique order-isomorphism between $\mathsf{ran}(g')$ and an ordinal (clearly $g''$ is a bijective $\in$-homorphism from $(\mathsf{TC}(q')\setminus\mathsf{dom}(f_i))$ to an ordinal),
    \item $f_{i+1}$ expands $f_i$ by putting $f_{i+1}=\alpha_i+g''(x)$ for $x\in\mathsf{dom}(g'')$ and then putting $f_{i+1}(q')=\alpha_i+\mathsf{ran}(g'')$.
    \end{itemize}
\end{itemize}
Finally we put $\mathsf{oemb}(e,s,q)={\bigcup_{i<\omega}} f_i$.

Now we are ready to compute $t^c(e,s)$. Let $f=\mathsf{oemb}(e,s,s)$ and $\alpha=\mathsf{rng}(f)$. We define the sequence of expanding partial enumerating functions $(e_\beta)_{\beta< \alpha}$ by recursion on $\beta< \alpha$. We put $e_{\beta}=u^c(e_{<\beta},\mathsf{TC}(\{\{t^p(e_{<\beta},y',\vec{z})\mid y'\in f^{-1}(\beta)\text{ and }\vec{z}\in s \}\}\cup s)),$ where $e_{<\beta}=e\cup \bigcup_{\gamma<\beta} e_{\gamma}$. Notice that by construction $e_{<\beta}$ is such that $t^d(e_{<\beta},y,\vec{z})=1$ for all $\vec{z}\in s$ and $y\in s$ with $f(y)<\beta$.  We output $e_{<\alpha}$ as the value of $t^c(e,s)$.
\end{proof}

For all $\Delta_0$-formulas $\varphi(\vec{x})$ of the language of $\prsou(\mathsf{enu})$ we could naturally define $\prsou$-terms $\varphi^d(e,\vec{x})$ and $\varphi^p(e,\vec{x})$ such that they are always equal to either $0$ or $1$, $\varphi^d(e,\vec{x})=1$ iff $e$ has enough information to evaluate the truth value of $\varphi(\vec{x})$, and if $\varphi^d(e,\vec{x})=1$, then $\varphi^p(e,\vec{x})$ computes this truth value. Formally we make the definition by recursion on construction of $\Delta_0$-formulas as follows:
\begin{itemize}
    \item If $\varphi(\vec{x})$ is $t(\vec{x})=u(\vec{x})$ or $t(\vec{x})\in u(\vec{x})$, then $\varphi^d(e,\vec{x})=t^d(e,\vec{x})u^d(e,\vec{x})$.
    \item If $\varphi(\vec{x})$ is $t(\vec{x})=u(\vec{x})$, then \[\varphi^p(e,\vec{x})=\begin{cases} 1& \text{if $t^p(e,\vec{x})=u^p(e,\vec{x})$}\\ 0 & \text{otherwise.}\end{cases}\]
    \item If $\varphi(\vec{x})$ is $t(\vec{x})\in u(\vec{x})$, then \[\varphi^p(e,\vec{x})=\begin{cases} 1& \text{if $t^p(e,\vec{x})\in u^p(e,\vec{x})$}\\ 0 & \text{otherwise.}\end{cases}\]
    \item If $\varphi(\vec{x})$ is $\lnot \psi(\vec{x})$, then $\varphi^d(e,\vec{x})=\psi^d(e,\vec{x})$ and $\varphi^p(e,\vec{x})=1-\psi^p(e,\vec{x})$.
    \item If $\varphi(\vec{x})$ is $\psi_1(\vec{x})\land \psi_2(\vec{x})$, then $\varphi^d(e,\vec{x})=\psi_1^d(e,\vec{x})\psi_2^d(e,\vec{x})$ and $\varphi^p(e,\vec{x})=\psi_1^p(e,\vec{x})\psi_2^p(e,\vec{x})$.
    \item If $\varphi(\vec{x})$ is $\forall y\in t(\vec{x})\;\psi(\vec{x},y)$, then \[\varphi^d(e,\vec{x})=t^d(e,\vec{x})\prod_{y\in t^p(e,\vec{x})}\psi^d(e,\vec{x},y)\] and \[\varphi^p(e,\vec{x})=\prod_{y\in t^p(e,\vec{x})}\psi^p(e,\vec{x},y).\]    
\end{itemize}

\begin{Lemma}
    For any $\Delta_0$-formula $\varphi(\vec{x})$ of the language of $\prsou(\mathsf{enu})$ it is $\prsou$-provable that for any partial enumerating function $e$, if $\varphi^d(e,\vec{x})=1$, then for any partial enumerating function $e'\supseteq e$, $e'\Vdash \varphi(\vec{x})$ iff $\varphi^p(e,\vec{x})=1$.
\end{Lemma}
\begin{proof}
    By straightforward induction on construction of $\Delta_0$-formulas.
\end{proof}

\begin{Lemma} For any $\Delta_0$-formula $\varphi(\vec{x})$ of the language of $\prsou(\mathsf{enu})$ there is a $\prsou$-term $\varphi^c(e,s)$ such that $\prsou$ proves that if $e$ is a partial enumerating function, $s\in\mathsf{dom}(e)$, $s$ is transitive, and $\vec{x}\in s$, then $\varphi^d(\varphi^c(e,s),\vec{x})=1$ and $\varphi^c(e,s)\supseteq e$.\end{Lemma}
\begin{proof}
We reason by induction on construction of $\Delta_0$-formulas.

If $\varphi(\vec{x})$ is either $t(\vec{x})\in u(\vec{x})$ or $t(\vec{x})= u(\vec{x})$, then $\varphi^c(e,s)=u^c(t^c(e,s),s)$.

If $\varphi(\vec{x})$ is $\lnot\psi(\vec{x})$, then $\varphi^c(s)=\psi^c(s)$.

If $\varphi(\vec{x})$ is $\psi_1(\vec{x})\land\psi_2(\vec{x})$, then $\varphi^c(e,s)=\psi_2(\psi_1^c(e,s),s)$.

If $\varphi(\vec{x})$ is $\forall y\in t(\vec{x})\,\psi(\vec{x})$, then \[\varphi^c(e,s)=\psi^c(t^c(e,s),\mathsf{TC}(s\cup\{t^p(t^c(e,s),\vec{x})\mid \vec{x}\in s\})).\]
\end{proof}

Now we are ready to prove Lemma \ref{forcing_axioms}:
\begin{proof}
    We only need to verify that the following three groups of axioms are forced: 
    \begin{enumerate}
        \item $\Delta_0$-foundation:\[\forall y (\forall z\in y\;\varphi(z)\to\varphi(y))\to\varphi(x)\text{, for $\Delta_0$ formulas $\varphi$.}\]
        \item Defining axiom for the function constructed by primitive set recursion:
        \begin{align*}
        \exists w\Big (\big ( \forall v\in w\exists u\in y\; v=R[t(x,y,\vec{z})](u,\vec{z})\big)\\\land\;\big( \forall u\in y\exists v\in w\;v=R[t(x,y,\vec{z})](u,\vec{z})\big )\\\land R[t(x,y,\vec{z})](y,\vec{z})=t(w,y,\vec{z})\Big ).    
        \end{align*}
        \item The defining axiom for $\mathsf{enu}$:
        \begin{equation}
        \label{def_enu_inst}
         x\ne \emptyset\land x\not\in U \to \forall x'\in x\exists y\in\omega(\mathsf{enu}(x,y)=x')\land \forall y\in \omega(\mathsf{enu}(x,y)\in x).   
        \end{equation}
    \end{enumerate}

    Note that $\Delta_0$-foundation schemata is implied over first-order logic by the following scheme:
    \begin{equation}\label{foundation_alt}
     \forall y\in x (\forall z\in y\;\varphi(z)\to\varphi(y))  \to \forall y\in x\, \varphi(y)\text{, for $\Delta_0$ formulas $\varphi$.}
    \end{equation}
    
    We reason in $\prsou+\text{``every set is countable''}$ to verify that $\emptyset$ forces (\ref{foundation_alt}) (for particular $\varphi$, $x$ and other parameters $\vec{p}$ of $\varphi$). For this it is enough to prove that any partial enumerating function $e$ could be extended to a partial enumerating  function $e'$ that forces (\ref{foundation_alt}). Let $s=\mathsf{TC}(\{x,\vec{p}\})$. Using an enumeration of $s$ we construct a partial enumerating function $e''$ covering $s$. We put $e'=\varphi^c(e'',s)$.  
    
    Now for any $r\supseteq e'$ and any $y\in x$ we have:
    \begin{equation} \label{deciding_extension} e'\Vdash \varphi(y)\iff r\Vdash \varphi(y)\iff \varphi^p(e',y,\vec{p})=1.\end{equation} 
    Using the first equivalence from (\ref{deciding_extension}) we observe that $e'$ forces (\ref{foundation_alt}) iff
    \[\forall y\in x (\forall z\in y\;e'\Vdash\varphi(z)\to e'\Vdash\varphi(y))  \to \forall y\in x\, e'\Vdash\varphi(y).\]
    And by the second equivalence from (\ref{deciding_extension}) we see that the assertion that $e'$ forces (\ref{foundation_alt}) is equivalent to
    \[\forall y\in x (\forall z\in y\;\varphi^p(e',z,\vec{p})=1\to \varphi^p(e',y,\vec{p})=1)  \to \forall y\in x\, \varphi^p(e',y,\vec{p})=1,\]
    which is clearly provable using $\Delta_0$-foundation.

    Our next goal is to show in $\prsou+\text{``every set is countable''}$ that $\emptyset$ forces some fixed instance of the defining axiom for the functions defined by primitive set recursion. Similarly to the verification fo $\Delta_0$-foundation we will need to see that for any partial enumerating function $e$ there is a partial enumerating function $e'\supseteq e$ that forces this instance. We put $e'=(R[t(x,y,\vec{z})])^c(e'',\mathsf{TC}(\{y,\vec{z}\}))$, where $e''$ is an extension of $e$ that covers all sets in $\mathsf{TC}(\{y,\vec{z}\}$ and the transitive closure itself. We observe that for this choice of $e'$ we have $(R[t(x,y,\vec{z})])^d(e',u,\vec{z})=1$ for all $u\in y\cup \{y\}$ and thus for any $r$ extending $e'$ and $u\in y$ we have 
    \[r\Vdash v=R[t(x,y,\vec{z})](u,\vec{z})\iff v=(R[t(x,y,\vec{z})])^p(e',u,\vec{z}).\]
    Hence to validate the forceability of (\ref{foundation_alt}) we pick $w=\{(R[t(x,y,\vec{z})])^p(e,u,\vec{z})\mid u\in y\}$. Then using the equivalence above we deduce that $e'\Vdash \forall v\in w\exists u\in y\; v=R[t(x,y,\vec{z})](u,\vec{z})$ and  $e'\Vdash \forall u\in y\exists v\in w\; v=R[t(x,y,\vec{z})](u,\vec{z})$. Finally, by definition of $(R[t(x,y,\vec{z})])^d$ the fact that $(R[t(x,y,\vec{z})])^d(e',y,\vec{z})=1$ implies that $t^d(e',w,y,\vec{z})=1$, which as before allows us to conclude that $e'\Vdash R[t(x,y,\vec{z})](y,\vec{z})=t(w,y,\vec{z})$, which finishes the verification of defining axiom for functions defined by primitive set recursion.

    Finally we defining axiom for $\mathsf{enu}$ in the same way. Now we simply pick $e'$ to be extension of the starting partial enumerating function $e$ by an enumeration of $x$. It is straightforward to check that this $e'$ forces (\ref{def_enu_inst}).
\end{proof}

\bibliography{References}

\end{document}